\documentclass{IEEEtran}

\usepackage{authblk}
\usepackage[round]{natbib}
\usepackage{hyperref}
\usepackage{url}
\usepackage[ruled,vlined,algo2e]{algorithm2e}
\usepackage{mathtools}
\usepackage{amsmath}

\usepackage{xcolor}
\hypersetup{
    colorlinks = true,
    linkcolor =  cyan,
    linkbordercolor = {white},
    citecolor=cyan,
}
\usepackage{amssymb}    % if you want extra symbols
\usepackage{mathrsfs}
\usepackage{cleveref}
\usepackage{empheq}
\usepackage{graphicx} % more modern
\usepackage{caption}
\usepackage{subcaption} 
\usepackage{multirow}
\newcommand{\ie}[0]{\emph{i.e.},~}
\newcommand{\eg}[0]{\emph{e.g.},~}

\newcommand{\defeq}{\ensuremath{\overset{\text{def}}{=}}}

\newcommand{\PP}[0]{\ensuremath{{p}}}%variable x
\newcommand{\msg}[0]{\ensuremath{\mathfrak{m}}}%variable x
\newcommand{\msgt}[2]{\ensuremath{\mathfrak{m}_{#1 \to #2}}}%variable x
\newcommand{\PPX}[0]{\ensuremath{{\PP^{\XX}}}}%variable x
\newcommand{\PPO}[0]{\ensuremath{{\PP^{\OO}}}}%variable x
\newcommand{\PPY}[0]{\ensuremath{{\PP^{\YY}}}}%variable x
\newcommand{\dd}[0]{\ensuremath{{d}}}
\newcommand{\matrixtype}[1]{\ensuremath{ {{\uppercase{{#1}}}}}}
\newcommand{\elementtype}[1]{\ensuremath{{\uppercase{{#1}}}}}
\newcommand{\XX}[0]{\ensuremath{\matrixtype{X}}}
\newcommand{\WW}[0]{\ensuremath{\matrixtype{W}}}
\newcommand{\YY}[0]{\ensuremath{\matrixtype{Y}}}
\newcommand{\ZZ}[0]{\ensuremath{\matrixtype{Z}}}
\newcommand{\OO}[0]{\ensuremath{\matrixtype{O}}}
\newcommand{\QQ}[0]{\ensuremath{\matrixtype{Q}}}
\newcommand{\OOO}[0]{\ensuremath{\mathcal{O}}}
\newcommand{\UU}[0]{\ensuremath{\matrixtype{U}}}
\newcommand{\xX}[0]{\ensuremath{\elementtype{X}}}
\newcommand{\wW}[0]{\ensuremath{\elementtype{W}}}
\newcommand{\yY}[0]{\ensuremath{\elementtype{Y}}}
\newcommand{\zZ}[0]{\ensuremath{\elementtype{Z}}}
\newcommand{\qQ}[0]{\ensuremath{\elementtype{Q}}}
\newcommand{\oO}[0]{\ensuremath{\elementtype{O}}}

\newcommand{\indextype}[1]{\ensuremath{{#1}}}
\newcommand{\indextypec}[1]{\ensuremath{{#1}}}
\newcommand{\m}[0]{\ensuremath{\indextype{m}}}
\newcommand{\n}[0]{\ensuremath{\indextype{n}}}
\renewcommand{\k}[0]{\ensuremath{\indextype{k}}}
\newcommand{\M}[0]{\ensuremath{\indextypec{M}}}
\newcommand{\N}[0]{\ensuremath{\indextypec{N}}}
\newcommand{\K}[0]{\ensuremath{\indextypec{K}}}
\newcommand{\truemath}[0]{\ensuremath{\mathrm{true}}}
\newcommand{\falsemath}[0]{\ensuremath{\mathrm{false}}}
\newcommand{\nullmath}[0]{\ensuremath{\mathrm{null}}}
\newcommand{\booltimes}[0]{\ensuremath{\bullet}}
\newcommand{\xortimes}[0]{\ensuremath{\ast}}

\newcommand{\margX}[0]{\ensuremath{\Xi}}
\newcommand{\margY}[0]{\ensuremath{\Upsilon}}
\newcommand{\aA}[0]{\ensuremath{\Phi}}
\newcommand{\bB}[0]{\ensuremath{\Psi}}
\newcommand{\cCh}[0]{\ensuremath{\Gamma}}
\newcommand{\cC}[0]{\ensuremath{\hat{\cCh}}}
\newcommand{\aAh}[0]{\ensuremath{\hat{\Phi}}}
\newcommand{\bBh}[0]{\ensuremath{\hat{\Psi}}}

\newcommand{\magn}[1]{\textbf{#1}}
\newcommand{\damping}[0]{\ensuremath{{\lambda}}}
\newcommand{\ttt}[1]{\ensuremath{^{(\mathrm{#1})}}}
\newcommand{\lpx}[0]{\ensuremath{{\phi}}}
\newcommand{\maxz}[1]{\ensuremath{\big(#1\big)_{+}}}
\newcommand{\ident}[0]{\ensuremath{\mathbb{I}}}
\newcommand{\ftype}[1]{\ensuremath{\mathsf{#1}}}
\newcommand{\f}[0]{\ensuremath{{\ftype{f}}}}
\newcommand{\h}[0]{\ensuremath{{\ftype{h}}}}
\newcommand{\g}[0]{\ensuremath{{\ftype{g}}}}

\begin{document}

%%%%%%%%%%%%%%%%%%%%%%%%%%%%%%

%% For titles, only capitalize the first letter
%% \title{Almost sharp fronts for the surface quasi-geostrophic equation}

\title{Boolean Matrix Factorization and Noisy Completion via Message Passing}

\renewcommand\Affilfont{\itshape\small}
\author{Siamak Ravanbakhsh$^{*}$}
\author{Barnab\'{a}s P\'{o}czos$^{*}$}
\author{Russell Greiner$^{**}$}
\affil{(*) Carnegie Mellon University, Pittsburgh, PA, 15213}
\affil{(**) University of Alberta, Edmonton, AB, Canada}

\maketitle

\begin{abstract}
Boolean matrix factorization and Boolean matrix completion from noisy observations are 
desirable unsupervised data-analysis methods due to their interpretability,
 but hard to perform due to their NP-hardness.
We treat these problems
as maximum a posteriori inference problems in a graphical model and present a message passing 
approach that scales linearly with the number of observations and factors. 
Our empirical study demonstrates that message passing is able to recover 
low-rank Boolean matrices, in the boundaries of theoretically possible recovery and
compares favorably with state-of-the-art in real-world applications, 
such collaborative filtering with large-scale Boolean data.
\end{abstract}

%\keywords{| Matrix Factorization | Matrix Completion | Boolean | Belief Propagation | Message Passing | Collaborative Filtering}

\section*{}

A body of problems in machine learning, communication theory and combinatorial optimization 
%-- from variations of principal component analysis to linear codes and compressed sensing -- 
involve the product form $\ZZ = \XX \odot \YY$ where $\odot$ operation corresponds to a type of matrix multiplication and 
\begin{align*}
%\label{eq:z}
 \ZZ = \{\zZ_{\m,\n}\}^{\M \times \N},
\XX = \{\xX_{\m,\k}\}^{\M \times \K}, \YY = \{\yY_{\k,\n}\}^{\K \times \N}.
\end{align*}
Here, often one or two components out of three are (partially) known and the task is to recover the unknown component(s).

% The diversity of these problems comes from their difference in \magn{a)} the underlying algebra: the definition of $\odot$ and its domain; \magn{b)} missing and observed values: which one of the matrices $\XX$, $\YY$ and $\ZZ$ are observed, and whether there are missing values; \magn{c)} priors: $\PPX(\XX)$, $\PPY(\YY)$ and $\PPZ(\ZZ)$ and \magn{d)} relative value of $\M$, $\N$ and $\K$.
% In many of these problems, message passing techniques -- \ie variants of belief propagation -- are state-of-the-art, sometimes reaching the boundaries of theoretically possible. 

A subset of these problems, which
are most closely related to Boolean matrix factorization and matrix completion,
can be 
expressed over the Boolean domain -- \ie $\zZ_{\m,\n},\ \xX_{\m,\k}, \ \yY_{\k,\n}\ \in\ \{\falsemath, \,\truemath\} \cong \{0,\,1\}$.
The two most common Boolean matrix products used in such applications are
\begin{subequations}
  \label{eq:prod}
\begin{empheq}{align}
&\ZZ = \XX \booltimes \YY  \Rightarrow  \zZ_{\m,\n} = \bigvee_{\k=1}^{\K} \xX_{\m,\k} \wedge \yY_{\k,\n} \hspace{-.1in}\label{eq:booltimes}\\
&\ZZ = \XX \xortimes \YY  \Rightarrow  \zZ_{\m,\n} \equiv \big ( \sum_{\k=1}^{\K} \xX_{\m,\k} \wedge \yY_{\k,\n}\big)\hspace{-.1in} \mod 2 \hspace{-.1in}\label{eq:xortimes}
\end{empheq}
\end{subequations}
where we refer to \Cref{eq:booltimes} simply as \textit{Boolean product} and we distinguish \Cref{eq:xortimes} as \textit{exclusive-OR (XOR) Boolean product}. One may think of Boolean product as ordinary matrix product where the values that are larger than zero in the product matrix are set to one. 
Alternatively, in XOR product, the odd (even) numbers are identically set to one (zero) in the product matrix. 

%When $\ZZ$ and $\YY$ are column vectors (\ie $\N = 1$), and $\zZ_{\m,\n} = \truemath$ for all $\m,\n$, 
%the problem of finding $\YY$ from a given $\XX$
%corresponds to
%\textit{Boolean satisfiability} (\aka SAT); here message passing has been able to solve hard random SAT instances close to 
%the unsatisfiability transition~\citep{braunstein2005survey,ravanbakhsh2014perturbed}. 
%Using the XOR product of \Cref{eq:xortimes}, the problem of recovering $\YY$ becomes exclusive-OR satisfiability (\aka XOR-SAT).

This model can represent \textbf{Low Density Parity Check (LDPC)} coding using the XOR product, with $N = 1$. In LDPC, the objective is to transmit the data vector $\YY \in \{0,1\}^{\K}$ though a noisy channel. For this, it is encoded by \Cref{eq:xortimes}, where $\XX \in \{0,1\}^{\m \times \k}$ is the \textit{parity check matrix} and vector $\ZZ \{0,1\}^{\M}$ is then sent though the channel with a noise model $\PPO(\OO \mid \ZZ)$, producing observation $\OO$.
Message passing decoding has been able to transmit $\ZZ$ and recover $\YY$ from $\OO$ 
at rates close to the theoretical capacity of the communication channel~\citep{gallager1962low}.

LDPC codes are in turn closely related to the \textbf{compressed sensing}~\citep{donoho2006compressed} -- so much so that successful binary LDPC codes (\ie matrix $\XX$) have been reused for compressed sensing~\citep{dimakis2012ldpc}.
In this setting,  the column-vector $\YY$ is known to be $\ell$-sparse (\ie $\ell$ non-zero values)
which makes it possible to use \textit{approximate message passing}~\citep{donoho2009message} to recover 
$\YY$ using few noisy measurements $\OO$ -- that is $\M \ll \K$ and similar to LDPC codes, the \textit{measurement matrix} $\XX$ is known. 
When the underlying domain and algebra is Boolean (\ie \Cref{eq:booltimes}), the compressed sensing problem 
reduces to the problem of (noisy) \textbf{group testing}~\citep{du1993combinatorial}
\footnote{The intuition is that the non-zero elements of the vector $\YY$ identify the presence or absence of a rare property (\eg a rare disease or manufacturing defect), 
therefore $\YY$ is sparse. The objective is to find these non-zero elements (\ie recover $\YY$) by screening a few ($\M \ll \K$) ``subsets'' of elements of $\YY$.
Each of these $\YY$-bundles corresponds to a row of $\XX$ (in \Cref{eq:booltimes}).} 
where message passing has been successfully applied in this setting as well~\citep{atia2012boolean,sejdinovic2010note}.

These problems over Boolean domain are special instances of the problem of Boolean factor analysis in which $\ZZ$ is given, but not $\XX$ nor $\YY$.
%in which both $\XX$ and $\YY$ have a matrix form and neither is given.
Here, inspired by the success of message passing techniques in closely related problems over ``real'' domain, we derive message passing solutions to a novel graphical model for ``Boolean'' factorization and matrix completion,
and show that simple application of Belief Propagation~\citep[BP;][]{pearl1982reverend} to this graphical model favorably compares with the state-of-the-art
in both Boolean factorization and completion.

% The intuition is that the non-zero elements of the vector $\YY$ identify the presence or absence of a rare property (\eg a rare disease or manufacturing defect), 
% therefore $\YY$ is sparse. The objective is to find these non-zero elements (\ie recover $\YY$) by screening a few ($\M \ll \K$) ``subsets'' of elements of $\YY$.
% %\RG{Are you saying that the Y-features are bundled togethers?  What are elements of Y}{}
% Each of these $\YY$-bundles corresponds to a row of $\XX$ (in \Cref{eq:booltimes}), and 
% message passing has been successfully applied in this setting as well~\citep{atia2012boolean,sejdinovic2010note}. 

In the following, we briefly introduce the Boolean factorization and completion problems 
in \Cref{sec:factor_analysis} and \Cref{sec:applications} reviews the related work.
\Cref{sec:formulation} formulates both of these problems in a Bayesian framework using
a graphical model. The ensuing message passing solution is introduced in \Cref{sec:message_passing}. Experimental study of \Cref{sec:experiments} demonstrates
 that message passing is an efficient and effective method for performing Boolean matrix factorization and noisy completion.

%With this background we elaborate the  setting that we study in this paper.
% These problems over Boolean domain are special instances of the problem of Boolean factor analysis in which both $\XX$ and $\YY$ have a matrix form and neither is given.
% This paper presents an efficient message passing solution to this general setting.

\subsection{Boolean Factor Analysis}\label{sec:factor_analysis}
The umbrella term ``factor analysis'' refers to the unsupervised methodology of expressing a set of observations
in terms of unobserved factors~\citep{mcdonald2014factor}.%
\footnote{While some definitions restrict factor analysis to variables over continuous domain or even probabilistic models with Gaussian priors, we take a more general view.} 
In contrast to LDPC and compressed sensing, in factor analysis, only (a partial and/or distorted version of) the matrix $\ZZ$ is observed,
and our task is then to find $\XX$ and $\YY$ whose product is close to $\ZZ$.
When the matrix $\ZZ$ is partially observed, a natural approach to Boolean matrix \textit{completion} is to find sparse and/or low-rank Boolean factors that would lead us to
missing elements of $\ZZ$. 
%Here, finding the factors (\ie factorization) is constrained by 
%observed elements of $\ZZ$ and low-rank/sparsity constraints ensure non-degeneracy of
%the reconstruction. 
In the following we focus on the Boolean product of \Cref{eq:booltimes},
noting that message passing derivation for factorization and completion using the XOR product of \Cref{eq:xortimes} is similar.

The ``Boolean'' factor analysis -- including factorization and completion -- has a particularly appealing form. 
This is because the Boolean matrix $\ZZ$ is simply written as disjunction of Boolean matrices of rank one -- that is $\ZZ = \bigvee_{\k=1}^{\K} \XX_{:,k} \booltimes \YY_{k,:}$, where $\XX_{:,k}$ and $\YY_{k,:}$ are column vector and row vectors of $\XX$ and $\YY$ respectively. 
%indicate the presence of a particular factor (\aka basis \citep{miettinen2006discrete}/ concept / category~\citep{belohlavek2010discovery,keprt2004binary}) in observations, 
%where each factor is in turn a Boolean pattern. 

\subsubsection{Combinatorial Representation}
The combinatorial representation of Boolean factorization is the
\textbf{biclique cover} problem in a bipartite graph $\mathcal{G} = (\mathcal{A} \cup \mathcal{B}, \mathcal{E})$. 
Here a bipartite graph has two disjoint node sets $\mathcal{A}$ (s.t. $|\mathcal{A}| = \M$)
and $\mathcal{B}$ (s.t. $|\mathcal{B}| = \N$) where the only edges are between these two sets -- \ie  $\mathcal{E} \subseteq \{(a,b) \mid  a \in \mathcal{A},\; b \in \mathcal{B}\}$. 
In our notation $\ZZ \in \{0,1\}^{M \times N}$ represents the incident matrix of $\mathcal{G}$
and the objective of factorization is to cover (only) the edges using $\K$ bicliques (\ie \textit{complete} bipartite sub-graphs of $\mathcal{G}$). Here the $\k^{th}$ biclique is identified
with a subset of $\mathcal{A}$, corresponding to $\XX_{:,\k}$, the $\k^{th}$ column of $\XX$, and a subset of $\mathcal{B}$, $\YY_{k,:}$, corresponding to the $\k^{th}$ row of $\YY$ the Boolean product of which is a Boolean matrix of rank $1$. The disjunction of these rank 1 matrices is therefore a biclique covering of the incident matrix $Z$.

\section{Applications and Related Work}\label{sec:related}\label{sec:applications}
%\noindent \textbf{Boolean Factorization.}
Many applications of Boolean factorization are inspired by its formulation as \textbf{tiling problem}~\citep{stockmeyer1975set}.\footnote{Since rows and columns in the rank one Boolean product $\XX_{:,\k} \booltimes \YY_{:,\k}^{T}$ can be permuted to form a ``tile'' -- \ie a sub-matrix where all elements are equal and different from elements outside the sub-matrix -- the Boolean factorization can be seen as tiling of matrix $\ZZ$ with tiles of rank one.}  
Examples include mining of Boolean databases~\citep{geerts2004tiling} to role mining~\citep{vaidya2007role,lu2008optimal} to 
bi-clustering of gene expression data~\citep{zhang2010binary}.
Several of these applications are accompanied by a method for approximating the Boolean factorization problem.

The most notable of these is the \textbf{``binary'' factorization}\footnote{
Binary factorization is different from Boolean factorization in the sense that in contrast to Boolean factorization $1+1 \neq 1$. Therefore the factors $\XX$
and $\YY$ are further constrained to ensure that $\ZZ$ does not contain any values other than zeros and ones.} of \citet{zhang2010binary} that uses an alternating optimization 
method to repeatedly solve a penalized non-negative matrix factorization problem over real-domain, 
where the penalty parameters try to enforce the desired binary form.
Note that a binary matrix factorization is generally a more constrained problem than Boolean factorization and therefore
it also provides a valid Boolean factorization.%
%\footnote{Although, factor $\g$ in our graphical model (and the resulting message passing procedure) can change 
%to accommodate binary constraints, here we are focused on the ``Boolean'' factorization and completion.} 

Among the heuristics~\citep[\eg][]{keprt2004binary,belohlavek2007fast} that directly apply to Boolean factorization, 
the best known is the \textbf{Asso} algorithm of \citet{miettinen2006discrete}. Since Asso is incrmental in $\K$, it can efficiently 
use the Minimum Description Length principle to select the best rank $K$ by incrementing its value~\citep{miettinen2011model}.

%\noindent \textbf{Boolean Completion.}
An important application of Boolean matrix completion is in collaborative filtering with Boolean (\eg like/dislike) data,
where the large-scale and sparsely observed Boolean matrices in modern applications demands a scalable and accurate Boolean matrix completion method. 

One of the most scalable methods for this problem is obtained by modeling the problem as a \textbf{Generalized Low Rank Model}~\citep[GLRM;][]{udell2014generalized},
that uses proximal gradient for optimization.
Using logistic or hinge loss can enforce binary values for missing entries. Using
the \textit{hinge loss}, GLRM seeks
 $$\arg\min_{\XX, \YY}  \sum_{(\m,\n) \in \Omega} \big ( 1 - (\sum_{\k} \xX_{\m,\k} \yY_{\k,\n} ) ( 2 \oO_{\m ,\n} - 1 ) \big )_+$$,
where $(2 \oO_{\m ,\n} - 1)$ changes the domain of observations to $\{-1, +1\}$ and $\Omega$ is index-set of observed elements. 

In the \textbf{1-Bit matrix completion} of \citet{davenport20141}, the single bit observation $\oO_{\m,\n}$ from a hidden real-valued matrix $\QQ$ 
is obtained by sampling from a distribution with the cumulative distribution function $f(\qQ_{\m,\n})$ -- \eg $\f(\qQ_{\m,\n}) = (1 + \exp(-\qQ_{\m,\n}))^{-1}$. 
For application to Boolean completion, our desired Boolean matrix is $\ZZ = \mathbb{I}(f(\QQ) {\geq} .5)$.  
1-Bit completion then minimizes the likelihood of observed entries, while constraining the nuclear norm of $\QQ$ 
\begin{align}\label{eq:1bitcompletion}
\arg\min_{\QQ} \; &\sum_{(\m,\n) \in \Omega}  \bigg (\oO_{\m,\n} \log(f(\qQ_{\m,\n})) + \\
&\oO_{\m,\n} \log(1 - f(\qQ_{\m,\n})) \bigg ) \quad s.t. \;\| \QQ \|_{*} \leq \beta \sqrt{\K \M \N}, \notag
\end{align}
where $\beta > 0$ is a hyper-parameter.

In another recent work, \citet{maurus2014ternary} introduce a method of \textit{ternary} matrix factorization 
that can handle missing data in Boolean factorization through ternary logic. In this model, the ternary
matrix $\ZZ$ is factorized to ternary product of a binary matrix $\XX$ and a ternary basis matrix $\YY$.
%Although this method has interesting applications of its own, due to ternary factors it cannot be used for Boolean completion.
%Here, we introduce a simple method that is not only faster than all methods above -- with $\OO(\K \,\| \Omega\|)$ complexity -- 
%it performs very well in practice.

%the non-linearity of Boolean matrix product has 
%confined most techniques to heuristics and local search methods
% Despite these applications, the non-linearity of Boolean product has confined most tractable 
% techniques to greedy heuristics~\citep[\eg][]{miettinen2006discrete,keprt2004binary,belohlavek2007fast} and
%  (penalized) relaxation of the problem to the real domain~\citep{zhang2007binary}. 

\section{Bayesian Formulation}\label{sec:formulation}
%This section formulates both Boolean factorization and completion problems as maximum a posteriori (MAP) inference in 
%a probabilistic graphical model. 
Expressing factorization and completion problems as a MAP inference problem is not new~\citep[\eg][]{mnih2007probabilistic}, 
neither is using message passing as an inference technique for these problems~\citep{krzakala2013phase,parker2013bilinear,kabashima2014phase}. 
However, message passing has not been previously used to solve the ``Boolean'' factorization/completion problem.% that poses its own subtleties. 
%In particular by introducing auxiliary variables and simplifying BP messages we provide a simple message update, 
%whose cost (per iteration) % with the cost that
%is linear in the number of observed elements of the matrix and the \textit{Boolean rank}, $K$.% (which is also provided as an input).

%(\aka Schein rank), which is the
% smallest number $\K$ for which a decomposition of $\ZZ$ to $\XX \booltimes \YY$ is exact~\citep{kim1982boolean}. 
%For noisy observations, \citet{miettinen2011model} suggest using the minimum description length (MDL) principle to estimate $\K$ (also see \citet{tatti2006dimension} for alternative measures). 

%Here, we also consider ``approximate'' decompositions, and assume that the value of $\K$ is given as an input. 
To formalize approximate decompositions for Boolean data, we use a communication channel,
where we assume that the product matrix $\ZZ$ is communicated through a {noisy binary erasure channel}~\citep{cover2012elements} to produce the observation $\OO \in \{0,1,\nullmath\}^{\M \times \N}$ where $\oO_{\m,\n} = \nullmath$, means this entry was erased in the channel. 
This allows us to model matrix completion using the same formalism that we use for low-rank factorization. 

For simplicity, we assume that each element of $\ZZ$ is independently transmitted (that is erased, flipped or remains intact) through the channel, meaning the following conditional probability completely defines the \textbf{noise model}:
\begin{align}\label{eq:noise_model}
\PPO(\OO \mid \ZZ) = \prod_{\m,\n} \PPO_{\m,\n}(\oO_{\m,\n} \mid \zZ_{\m,\n})
\end{align}
Note that each of these conditional probabilities can be represented using six values -- one value per each pair of $\oO_{\m,\n} \in \{0,1,\nullmath\}$ and $\zZ_{\m,\n} \in \{0,1\}$. 
This setting allows \textit{the probability of erasure to depend on the value of $m$, $n$ and $\zZ_{\m,\n}$}. 

The objective is to recover $\XX$ and $\YY$ from $\OO$.
However, due to its degeneracy, recovering $\XX$ and $\YY$ is only up to a $\K \times \K$ permutation matrix $\UU$ -- that is
$\XX \booltimes \YY = (\XX \booltimes \UU) \booltimes (\UU^T \booltimes \YY)$. 
A Bayesian approach can resolve this ambiguity by defining non-symmetric \textbf{priors}
\begin{subequations}
  \label{eq:xypriors}
\begin{empheq}{align}
\PPX(\XX) =  \prod_{\m,\k}\PPX_{\m,\k}(\xX_{\m,\k}) \\ 
\PPY(\YY) = \prod_{\k,\n}\PPY_{\k,\n}(\yY_{\k,\n}) 
\end{empheq}
\end{subequations}
where we require the a separable product form for this prior. Using strong priors can enforce 
sparsity of $\XX$ and/or $\YY$, leading to well-defined factorization and completion problems where $K > M,N$.

Now, we can express the problem of recovering $\XX$ and $\YY$ as a 
\textit{maximum a posteriori} (MAP) inference problem \\
$\arg\max_{\XX, \YY}\; \PP(\XX, \YY \mid \OO)$, where the posterior is 
\begin{align}
\label{eq:map}
\PP(\XX, \YY \mid \OO) \;\propto\; \PPX(\XX) \; \PPY(\YY) \; \PPO(\OO \mid \XX \booltimes \YY)
%&\prod_{\m,\k}\PPX_{\m,\k}(\xX_{\m,\k})\; \prod_{\k,\n}\PPY_{\k,\n}(\yY_{\k,\n})\\ 
%&\prod_{\m,\n} \PPO_{\m,\n}\left(\oO_{\m,\n} \mid \bigvee_{\k=1}^{\K} \xX_{\m,\k} \wedge \yY_{\k,\n}\right)\notag
\end{align}

% where the marginal is
% \begin{align}
% \label{eq:marginal}
% \PP(\xX_{\m,\k}\mid \OO) = \sum_{\backslash \xX_{\m,\k}} \PP(\XX, \YY \mid \OO)
% \end{align}

% This formulation can represent
% \magn{(a)} low-rank Boolean matrix factorization, where $\K < \M,\N$;\  \
% \magn{(b)} Boolean \textit{dictionary learning}, where $\K > \M,\N$ and non-uniform priors encouraging sparsity of $\XX$ or $\YY$; 
% \magn{(c)} Boolean (noisy) \textit{matrix completion}, where only a fraction of the (noisy) product matrix are observed -- \ie $\PPO_{\m,\n}( 
% \oO_{\m,\n} =  % RG
% \nullmath \mid \zZ_{\m,\n}) > 0$. 

%Alternatively we may seek to recover $\XX$ and $\YY$ using the \textit{marginals} of the posterior 
%-- \eg use $\arg_{\xX_{\m, \k}}\max \;\PP(\xX_{\m,\k}\mid \OO)$ to recover $\xX_{\m,\k}$.

Finding the maximizing assignment for \Cref{eq:map} is ${NP}$-hard~\citep{stockmeyer1975set}. %Our treatment of SAT as a special case of factorization with $\m=1$ and fixed $\YY$ proves that even the problem of finding the Boolean matrix $\XX$ (when $\YY$ and $\ZZ$ are given) such that $\ZZ = \XX \booltimes \YY$,
%is $\mathbb{NP}$-hard. 
Here we introduce a graphical model to represent the posterior and
use a simplified form of BP to approximate the MAP assignment. 

An alternative to finding the MAP assignment is that of finding the marginal-MAP -- \ie $$\arg\max_{\XX_{\m,\k}} \PP(\xX_{\m,\k} \mid \OO) = \arg\max_{\XX_{\m,\n}} \sum_{\XX \backslash \xX_i, \YY} \PP(\XX, \YY \mid \OO).$$
While the MAP assignment is the optimal ``joint'' assignment to $\XX$ and $\YY$, finding the marginal-MAP corresponds to optimally estimating individual assignments for each variable, while the other variable assignments are marginalized. We also provide the message passing solution to this alternative in Appendix B. 
%Here, we focus on the MAP estimate as we find it more useful in practice.

\subsection{The Factor-Graph}\label{sec:factor_graph}
\begin{figure}
  \begin{center}
    \includegraphics[width=0.25\textwidth]{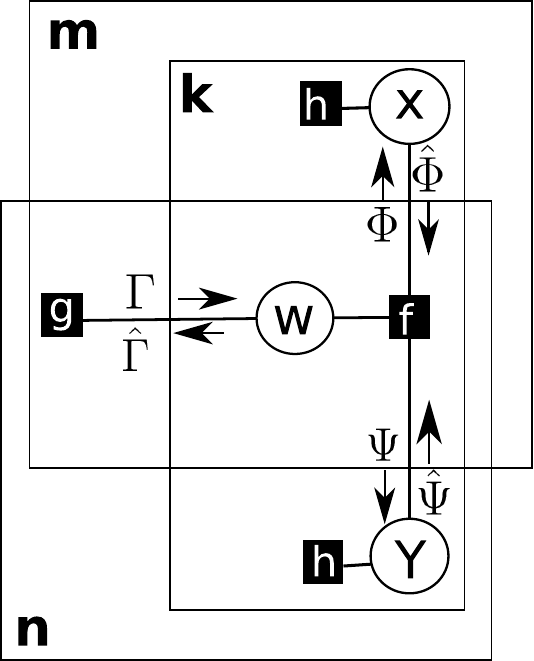}
  \end{center}
\vspace*{-1em}
  \caption{\it \small{The factor-graph and the message exchange between variables and factors.}}\label{fig:fg}
\vspace*{-1em}
\end{figure}

\Cref{fig:fg} shows the factor-graph \citep{kschischang2001factor} representation of the posterior \Cref{eq:map}. 
Here, variables are circles and factors are squares. The factor-graph is a bipartite graph, connecting each factor/function to its relevant variables. 
This factor-graph has one variable $\xX_{\m,\k} \in \{0,1\}$ for each element of $\XX$,  and a variable $\yY_{\k,\n} \in \{0,1\}$ for each element of $\YY$. 
In addition to these $\K \times (\M + \N)$ variables, we have introduced $\K \times \M \times \N$ auxiliary variables $\wW_{\m, \n, \k} \in \{0,1\}$. For Boolean matrix completion the number of auxiliary variables is $\K |\Omega|$, where $\Omega = \{(\m,\n) | \oO_{\m,\n} \neq \nullmath \}$ is the set of observed elements (see \Cref{sec:simplification}).

We use plate notation (often used with directed models) in representing this factor-graph.
\Cref{fig:fg} has three plates for $1 \leq \m \leq \M$, $1 \leq \n \leq \N$ and $1 \leq \k \leq \K$ (large transparent boxes in \Cref{fig:fg}).
In plate notation, all variables and factors on a plate are replicated. For example, variables on the $m$-plate are replicated for $1 \leq \m \leq \M$. 
Variables and factors located on more than one plate are replicated for all combinations of their plates.
For example, since variable $\xX$ is in common between $\m$-plate and $\k$-plate, it refers to $\M \times \K$ binary variables -- \ie $\xX_{\m,\k} \; \forall \m,\k$.  

\subsubsection{Variables and Factors}
The auxiliary variable $\wW_{\m,\n,\k}$  represents the Boolean product of $\xX_{m,\k}$ and $\yY_{\k,\n}$ -- 
\ie $\wW_{\m,\n,\k} = \xX_{m,\k} \wedge \yY_{\k,\n}$. This is achieved through $\M \times \N \times \K$ 
hard constraint factors
\begin{align*}
\f_{\m,\n,\k}(\xX_{\m,\k}, \yY_{\k,\n}, \wW_{\m,\n,\k}) \; = \; \ident(\wW_{\m,\n,\k} = \xX_{m,\k} \wedge \yY_{\k,\n}) 
\end{align*}
where $\ident(.)$ is the identity function on the inference semiring~\citep[see][]{ravanbakhsh2014revisiting}.
For the max-sum inference $\ident_{\text{max-sum}}(\truemath) = 0$ and $\ident_{\text{max-sum}}(\falsemath) = -\infty$. 
%Analogously, for sum-product inference
%$\ident_{\text{sum-prod}}(\truemath) = 1$ and $\ident_{\text{sum-prod}}(\falsemath) = 0$.

Local factors 
%\begin{align*}
$\h_{\m,\k}(\xX_{\m,\k}) \; = \; \log(\PPX(\xX_{\m,\k}))$ %\\%\qquad %\xX_{\m,\k} \, \xR_{\m,\n} \quad \text{and} \quad 
and
$\h_{\k,\n}(\yY_{\k, \n}) \; = \; \log(\PPY(\yY_{\k,\n}))$
%\yY_{\m,\k} \, \yR_{\m,\n}
%\end{align*}
represent the logarithm of priors over $\XX$ and $\YY$ in \Cref{eq:map}.

Finally, the noise model in \Cref{eq:map} is represented by $\M\times\N$ factors over auxiliary variables %-- \ie $\PP(\oO_{\m,\n} \mid \bigvee_{\k=1}^{\K} \xX_{\m,\k} \wedge \yY_{\k,\n})$ 
\begin{align*}
  \g_{\m,\n}(\{\wW_{\m,\n,\k}\}_{1\leq \k\leq \K})  = 
\log\bigg(\PPO_{\m,\n}(\oO_{\m,\n} \mid \bigvee_{\k} \wW_{\m,\n,\k})\bigg).% \; = \; \oR_{\m,\n}\; (-1)^{1 + (\bigvee_{\k} \wW_{\m,\n,\k})}
\end{align*}
%where $\sim x$ is the logical not of $x$ and we used $[\falsemath, \,\truemath] \cong [0,\,1]$ in our notation.
Although our introduction of auxiliary variables is essential in building our model, the factors of this type have been used in the past. %re-appeared almost independently in different application domains in the past.
In particular, factor $\g$ is generalized by a high-order family of factors with tractable inference, known as cardinality-based potentials~\citep{gupta2007efficient}.
This factor is also closely related to noisy-or models~\citep{pearl2014probabilistic, middleton1991probabilistic}; 
where MCMC~\citep{wood2012non} and variational inference~\citep{vsingliar2006noisy} has been used to solve more sophisticated probabilistic models of this nature.

The combination of the factors of type $\g$ and $\f$, represent the term 
$\PP(\oO_{\m,\n} \mid \bigvee_{\k=1}^{\K} \xX_{\m,\k} \wedge \yY_{\k,\n})$ in \Cref{eq:map} and the local factors $\h$, 
represent the logarithm of the priors. It is easy to see that the sum of all the factors above, evaluates to
 the logarithm of the posterior 
\begin{align*}
\log(\PP(\XX, \YY \mid \OO) &= \sum_{{\m,\k}} \h_{\m,\k}(\xX_{\m,\k}) + \sum_{\k,\n} \h_{\k, \n}(\xX_{\k, \n}) \\
&+ \sum_{\m,\n} \g_{\m,\n}(\{\xX_{m,\k} \wedge \yY_{\k,\n}\}_{1\leq \k\leq \K})
\end{align*}
 if  $\wW_{\m,\n,\k} = \xX_{m,\k} \wedge \yY_{\k,\n}\, \forall \m,\n,\k$ and $-\infty$ otherwise.
Therefore, maximizing the sum of these factors is equivalent to MAP inference for \Cref{eq:map}. 
% Note that expressing factorization and completion problems as a MAP inference problem is not new
% neither is using message passing
% to solve factorization problem~\citep[\eg]{kabashima2014phase}. However, subtleties such as introducing auxiliary variables and simplification
% in the following section is essential in approximating the MAP problem using BP.

%Once we have the factor-graph representing the posterior, we can use max-sum Belief Propagation~(BP) to approximate the MAP assignment. As \Cref{fig:fg} shows, the messages are passed between nodes and factors in both directions.
%and 
%\Cref{eq:map_update} and \Cref{eq:marg_update} are the simplified form of these BP message updates. 

\section{Message Update}\label{sec:algorithm}\label{sec:message_passing}
%This section introduces message passing solutions to approximate this MAP inference  (see \Cref{sec:marginal} for marginalization algorithm).
% Before writing the updates, we introduce some short notation for log-ratio of priors as well as the likelihood (\ie noise model probabilities):
% \begin{subequations}
% \begin{empheq}{align}
% \xR_{\m,\k} \quad &= \quad \log \bigg(\frac{\PPX_{\m,\k}(1)}{\PPX_{\m,\k}(0)}\bigg )\\
% \yR_{\k,\n} \quad &= \quad \log \bigg(\frac{\PPY_{\k,\n}(1)}{\PPY_{\k,\n} (0)}\bigg )\\
% \oR_{\m,\n} \quad &= \quad \log \bigg (\frac{\PPO_{\m,\n}(\oO_{\m,\n} \mid 1)}{\PPO_{\m,\n}(\oO_{\m,\n} \mid 0)} \bigg ) \label{eq:oRmn}
% \end{empheq}
% \end{subequations}
\SetAlgoNoLine
\SetInd{.05em}{.01em}
\begin{algorithm2e}[t!]
\caption{message passing for Boolean matrix factorization/completion}
\label{alg:1}
\DontPrintSemicolon
\small
\KwIn{1) observed matrix $\OO \in \{0,1\}^{\M \times \N}\;\forall m,n$;\\ 2) $\K \in \mathbb{N}$;\\ 3) priors  $\PPX_{\m,\k}, \PPY_{\n,\k}\;\forall m,n,k$;\\4) noise model $\PPO_{\m,\n} \;\forall m,n,k$ } %5) damping $\damping \in (0,1]$.}
\KwOut{$\XX \in \{0,1\}^{\M \times \K}$ and $\YY \in \{0,1\}^{\K \times \N}$.}
$t := 0$\\
\textbf{init} $\aA_{\m,\n,\k}\ttt{t},\bB_{\m,\n,\k}\ttt{t},\aAh_{\m,\n,\k}\ttt{t},\bBh_{\m,\n,\k}\ttt{t}, \cC_{\m,\n,\k}\ttt{t}$ and $\cCh_{\m,\n,\k}\ttt{t} \;\;\;\; \forall \m,\n,\k$ \\% 1\leq\m\leq \M,\,1 \leq \n \leq \N,\, 
%1 \leq \k \leq \K$.\\ 
%using $\log(U) - \log(1-U),\; U \sim \mathrm{Uniform}(0,1)$\\
%\begin{algorithmic}
\While{$t < T_{\max}$ \textbf{and} not converged for all $m,n,k$}{
\begin{subequations}
\label{eq:map_update}
\begin{empheq}{align}
%\begin{align}
&\aA_{\m,\n,\k}\ttt{t+1}  :=  \maxz{\cCh_{\m,\n,\k}\ttt{t} + \bBh_{\m,\n,\k}\ttt{t}} \, - \, \maxz{\bBh_{\m,\n,\k}\ttt{t}} \label{eq:xin}\\
&\bB_{\m,\n,\k}\ttt{t+1} :=  \maxz{\cCh_{\m,\n,\k}\ttt{t} + \aAh_{\m,\n,\k}\ttt{t}} \, - \, \maxz{\aAh_{\m,\n,\k}\ttt{t}} \label{eq:yin}\\
&\aAh_{\m,\n,\k}\ttt{t+1} :=  \log \bigg (\frac{\PPX_{\m,\k}(1)}{\PPX_{\m,\k}(0)} \bigg) + \sum_{\n' \neq \n} \; \aA_{\m,\n',\k}\ttt{t}  \label{eq:xout}\\
&\bBh_{\m,\n,\k}\ttt{t+1} :=   \log \bigg(\frac{\PPY_{\n,\k}(1)}{\PPY_{\n,\k}(0)} \bigg ) + \sum_{\m' \neq \m} \; \bB_{\m',\n,\k}\ttt{t}  \label{eq:yout}\\
&\cC_{\m,\n,\k}\ttt{t+1} := \min \bigg \{ \aAh_{\m,\n,\k}\ttt{t} + \bBh_{\m,\n,\k}\ttt{t}, \notag \\ 
&\quad\quad\quad\quad\quad\quad\quad \aAh_{\m,\n,\k}\ttt{t}, \bBh_{\m,\n,\k}\ttt{t} \bigg\}\hspace{-.1in}\label{eq:zout} \\
&\cCh_{\m,\n,\k}\ttt{t+1} := \min \bigg \{ \maxz{-\max_{\k' \neq \k}\cC_{\m,\n,\k'}\ttt{t}}, \notag \\
&\sum_{\k' \neq \k} \maxz{\cC_{\m,\n,\k'}\ttt{t}} + \log \bigg(\frac{\PPO_{\m,\n}(\oO_{\m,\n} \mid 1)}{\PPO_{\m,\n}(\oO_{\m,\n} \mid 0)} \bigg ) \bigg\}\hspace{-.1in}\label{eq:zin}
\end{empheq}
\end{subequations}
%\end{align}
}

\textbf{calculate} log-ratio of the posterior marginals\\
\begin{subequations}
\label{eq:marg}
\begin{empheq}{align}
  \margX_{\m,\k} := \log\bigg(\frac{\PPX_{\m,\k}(1)}{\PPX_{\m,\k}(0)}\bigg) + \sum_{\n} \aA_{\m,\n,\k}\ttt{t} \label{eq:margx}\\
  \margY_{\k, \n} := \log\bigg (\frac{\PPY_{\k, \n}(1)}{\PPY_{\k, \n}(0)}\bigg) + \sum_{\m} \bB_{\m,\n,\k}\ttt{t} \label{eq:margy}
\end{empheq}
\end{subequations}

\textbf{calculate} $\XX$ and $\YY$\\
\begin{subequations}
\label{eq:threshold}
\begin{empheq}{align}
  &\xX_{\m,\k} :=  \begin{cases}
1, & \text{if} \;  \margX_{\m,\k}  > 0 \label{eq:margx_threshold}\\
0, & \text{otherwise}
\end{cases}\\
 &\yY_{\k, \n} :=  \begin{cases}
1, & \text{if} \;  \margY_{\k, \n}  > 0 \label{eq:margy_threshold}\\
0, & \text{otherwise}
\end{cases}
\end{empheq}
\end{subequations}

\Return{\XX,\YY}
\end{algorithm2e}

Max-sum Belief Propagation (BP) is a message passing procedure for approximating the MAP assignment in a graphical model.
In factor-graphs without loops, max-sum BP is simply an exact dynamic programming approach that leverages the distributive law.
In loopy factor-graphs the approximations of this message passing procedure 
is justified by the fact that it represents the zero temperature limit to the sum-product BP, which is in turn a fixed point iteration procedure whose
fixed points are the local optima of the Bethe approximation to the free energy~\citep{yedidia2000generalized}; see also \citep{weiss2012map}.
For general factor-graphs, it is known that the approximate MAP solution obtained using max-sum BP is optimal within its ``neighborhood''~\citep{weiss2001optimality}.
%Although few applications of max-sum BP to loopy factor-graphs have optimality guarantees, approximations using BP have 
%been successfully applied to machine learning, coding theory and 
%combinatorial optimization problems. % in \citep[\eg][]{frey2007clustering}

We apply max-sum BP to approximate the MAP assignment of the factor-graph of \Cref{fig:fg}.
This factor-graph is very densely connected and therefore, one expects BP to oscillate or fail to find a good solution.
However, we report in \Cref{sec:experiments} that BP performs surprisingly well. This can be attributed to the week influence of majority of the factors, 
often resulting in close-to-uniform messages.
Near-optimal behavior of max-sum BP in dense factor-graph is not without precedence \citep[\eg][]{frey2007clustering,decelle2011asymptotic,ravanbakhsh2014augmentative}.

The message passing for MAP inference of \Cref{eq:map} involves message exchange between all variables and their neighboring factors
in both directions. Here, each message is a Bernoulli distribution. For example $\msgt{\xX_{\m,\k}}{\f_{\m,\n,\k}}(\xX_{\m,\n}): \{0,1\} \to \Re^2$
is the message from variable node $\xX_{\m,\n}$ to the factor node $\f_{\m,\n,\k}$.  For binary variables, it is convenient to work with the log-ratio of messages -- \eg we use
$\aAh_{\m,\n,\k} = \log \big(\frac{\msgt{\xX_{\m,\k}}{\f_{\m,\n,\k}}(1)}{\msgt{\xX_{\m,\k}}{\f_{\m,\n,\k}}(0)}\big)$ and the log-ratio of the message is opposite direction is denoted by $\aAh$.
Messages $\bB$, $\bBh$, $\cC$ and $\cCh$ in \Cref{fig:fg} are defined similarly.
%three message types $\aA, \bB$ and $\cC$, 
%where each message is a real value indexed by $1 \leq \m \leq \M$, $1 \leq \n \leq \N$ and $1 \leq \k \leq \K$.
For a review of max-sum BP and the detailed derivation of the simplified BP updates for this factor-graph, see Appendix A. In particular, a naive application of BP to obtain 
messages $\cCh_{\m,\n}$ from the likelihood factors  $\g_{\m,\n}(\{\wW_{\m,\n,\k}\}_{1\leq \k\leq \K})\; \forall \m,\n$ to the auxiliary variables $\wW_{\m,\n,\k}$ has a $\mathcal{O}(2^{\K})$ cost. 
In Appendix A, we show how this can be reduced to $\mathcal{O}(\K)$.
\Cref{alg:1} summarizes the simplified message passing algorithm.

At the beginning of the {Algorithm}, $t = 0$, messages are initialized with some random value -- \eg using $\log(U) - \log(1 - U)$ where $U \sim \mathrm{Uniform(0,1)}$. 
Using the short notation 
$\maxz{a} \; = \; \max\{0, a\}$, at time $t+1$, the messages are updated using 1) the message values at the previous time step $t$; 2) the prior; 3) the noise model and observation $\OO$.
The message updates of \Cref{eq:map_update} are repeated until convergence or a maximum number of iterations $T_{\max}$ is reached.
We decide the convergence based on the maximum absolute change in one of the message types \eg $\max_{\m,\n,\k} |\aA_{\m,\n,\k}\ttt{t+1} - \aA_{\m,\n,\k}\ttt{t}| \overset{?}{\leq} \epsilon$.

Once the message update converges, at iteration $T$, we can use the values for $\aA_{\m,\n,\k}\ttt{T}$ and $\bB_{\m, \n,\k}\ttt{T}$ to recover 
the log-ratio of the marginals $\PP(\xX_{\m,\k})$ and $\PP(\yY_{\n,\k})$. These log-ratios are denoted by $\margX_{\m,\k}$ and $\margY_{\k,\n}$ in \Cref{eq:marg}.
A positive log-ratio $\margX_{\m,\k} > 0$ means $\PP(\xX_{\m,\k} = 1) > \PP(\xX_{\m,\k} = 0)$ and the posterior favors $\xX_{\m,\k} = 1$. In this way the marginals
are used to obtain an approximate MAP assignment to both $\XX$ and $\YY$.

For better convergence, we also use \textit{damping} in practice. For this, one type of messages is updated to a linear combination of messages at time $t$ and $t+1$
 using a \textit{damping parameter} $\damping \in (0,1]$. Choosing $\aAh$ and $\bBh$ for this purpose, the updates of \Cref{eq:xout,eq:yout} become
\begin{align}\label{eq:damping}
\aAh_{\m,\n,\k}\ttt{t+1} &:=  (1 - \damping) \aAh_{\m,\n,\k}\ttt{t} + \\
&\damping \bigg ( \log \bigg( \frac{\PPX_{\m,\k}(1)}{\PPX_{\m,\k}(0)} \bigg) + \sum_{\n' \neq \n}  \aA_{\m,\n',\k}\ttt{t} \bigg ),  \notag \\
\bBh_{\m,\n,\k}\ttt{t+1} &:= (1 - \damping) \bBh_{\m,\n,\k}\ttt{t}  + \notag \\
& \damping \bigg ( \log \bigg( \frac{\PPY_{\n,\k}(1)}{\PPY_{\n,\k}(0)} \bigg ) + \sum_{\m' \neq \m} 
\bB_{\m',\n,\k}\ttt{t} \bigg ). \notag
\end{align}

\subsection{Further Simplifications}\label{sec:simplification}
\magn{Partial knowledge.} If any of the priors, $\PP(\xX_{\m,\k})$ and $\PP(\yY_{\n,\k})$, are zero or one, it means that $\XX$ and $\YY$ are partially known. The message updates of \Cref{eq:xout,eq:yout} will assume $\pm \infty$ values, to reflect these hard constrains. In contrast, for uniform priors, the log-ratio terms disappear.

\magn{Matrix completion speed up.} Consider the case where 
$\log \big (\frac{\PPO(\oO_{\m,\n} \mid 1)}{\PPO(\oO_{\m,\n} \mid 0)} \big) = 0$ in \Cref{eq:zin}
--
\ie the probabilities in the nominator and denominator are equal.
An important case of this happens in matrix completion, when the probability of erasure is independent of the value of $\zZ_{\m,\n}$ -- 
that is $\PPO(\nullmath \mid \zZ_{\m,\n} = 0) = \PPO(\nullmath \mid \zZ_{\m,\n} = 1) = \PPO(\nullmath)$ for all $\m$ and $\n$. 

 It is easy to check that in such cases,
 $\cCh_{\m,\n,\k} = \min \big( \maxz{-\max_{\k' \neq \k}\cC_{\m,\n,\k}\ttt{t}} , \sum_{\k' \neq \k} \maxz{\cC_{\m,\n,\k}\ttt{t}} \big )$ is always zero. 
This further implies that $\aAh_{\m,\n,\k}$ and $\bBh_{\m,\n,\k}$ in \Cref{eq:xout,eq:yout} are also always zero and calculating $\cC_{\m,\n,\k}$ in \Cref{eq:zin} is pointless. 
The bottom-line is that we only need to keep track of messages where this log-ratio is non-zero. 
Recall that $\Omega =  \{(\m,\n) \mid \oO_{\m,\n} \neq \nullmath \}$ denote the observed entries of $\OO$. 
Then in the message passing updates of \Cref{eq:map_update} in \Cref{alg:1}, wherever the indices $\m$ and $\n$ appear, we may restrict them to the set $\Omega$.

% \magn{Noiseless channel.} 
% When $\PPO(1 \mid 1) = \PPO(0 \mid 0) = 1$,
% $\oR_{m,n}$ evaluates to $-\infty$ or $+\infty$ for $\oO_{m,n} = 0$ and $\oO_{m,n} = 1$ respectively. In the former case -- \ie $\oR_{m,n} = -\infty$ -- 
% the update of \Cref{eq:zin} is unnecessary as $\cCh_{\m,\n,\k}\ttt{t+1} = -\infty$ at all time. Consequently, \Cref{eq:xin,eq:yin} further simplify, speeding up the message update for the noiseless Boolean matrix factorization problem.
%especially when $\OO$ is sparse. 
%The practical implication is that for the noiseless case, we do not need to keep
%track of $\cC_{\m,\n,\k}$ and $\cCh_{\m,\n,\k}$ for $\$

\magn{Belief update.} 
Another trick to reduce the complexity of message updates is in calculating $\{\aAh_{\m,\n,\k}\}_{\n}$ and $\{\bBh_{\m,\n,\k}\}_{\m}$ in \Cref{eq:xout,eq:yout}. 
We may calculate the marginals $\margX_{\m,\k}$ and $\margY_{\k,\n}$ using \Cref{eq:marg},
and replace the \Cref{eq:damping}, the damped version of the \Cref{eq:xout,eq:yout}, with
\begin{subequations}
\label{eq:out_fast}
\begin{empheq}{align}
\aAh_{\m,\n,\k}\ttt{t+1}  &:=  (1 - \damping) \aAh_{\m,\n,\k}\ttt{t} + \damping \big (    
\margX_{\m,\k}\ttt{t} - \aA_{\m,\n,\k}\ttt{t} \big )\label{eq:xout_fast}\\
\bBh_{\m,\n,\k}\ttt{t+1}  &:=  (1 - \damping) \bBh_{\m,\n,\k}\ttt{t} + \damping \big ( \margY_{\k,\n}\ttt{t} - \bB_{\m,\n,\k}\ttt{t} \big )\label{eq:yout_fast}
\end{empheq}
\end{subequations}
where the summation over $\n'$ and $\m'$ in \Cref{eq:xout,eq:yout} respectively, is now performed only once (in producing the marginal) and reused.

\magn{Recycling of the max.} 
Finally, using one more computational trick the message passing cost is reduced to linear:
in \Cref{eq:zout}, 
% in the term $\maxz{-\max_{\k' \neq \k}\cC_{\m,\n,\k}\ttt{t}}$
the maximum 
of the term $\maxz{-\max_{\k' \neq \k}\cC_{\m,\n,\k}\ttt{t}}$
is calculated for each of $\K$ messages $\{\cCh_{\m,\n,\k}\}_{\k \in \{1,\ldots, \K\}}$. 
Here, we may calculate the ``two'' largest values in the set $\{\cC_{\m,\n,\k}\ttt{t}\}_{\k}$ only once and reuse them in the updated for all $\{\cCh_{\m,\n,\k}\}_\k$ -- 
\ie if the largest value is $\cC_{\m,\n,\k^*}\ttt{t}$ 
then 
we use the second largest value, only in producing $\cCh_{\m,\n,\k^*}$.

\magn{Computational Complexity.} 
All of the updates in (\ref{eq:xin},\ref{eq:yin},\ref{eq:zin},\ref{eq:zout},\ref{eq:out_fast}) have a constant computational cost. Since these are performed for $\K |\Omega|$ messages, and the updates in calculating the marginals \Cref{eq:margx,eq:margy} are $\OOO(\K |\Omega| )$, the complexity of one iteration is $\OOO(\K |\Omega|)$.

% Now let us consider the extreme cases of the noise model.  In a \textit{binary symmetric channel (BSC)}, all the elements of $\OO$ are observed -- \ie $\PP(\nullmath \mid \zZ_{\m,\n}) = 0$ for all $\m,\n, \zZ_{\m,\n} \in \{0,1\}$. Moreover, for BSC 
% $\log (\frac{\PP(1 \mid 1)}{\PP(1 \mid 0)}) = -\log (\frac{\PP(0 \mid 1)}{\PP(0 \mid 0)}) = c > 0$.
% Therefore, \Cref{eq:zin} simplifies to
% \begin{align}
% \cCh_{\m,\n,\k}\ttt{t+1} \quad &= \quad \min \bigg( 0, -\max_{\k' \neq \k}\cC_{\m,\n,\k}\ttt{t} , \sum_{\k' \neq \k} \maxz{\cC_{\m,\n,\k}\ttt{t}} + (2\oO_{\m,\n} -1) c \bigg)
% \end{align}

\begin{figure*}
  \centering
\hbox{
%\subcaptionbox*{}{\includegraphics[width=.3\textwidth]{figures/nimfa_rib.pdf}}
\subcaptionbox*{}{\includegraphics[width=\textwidth]{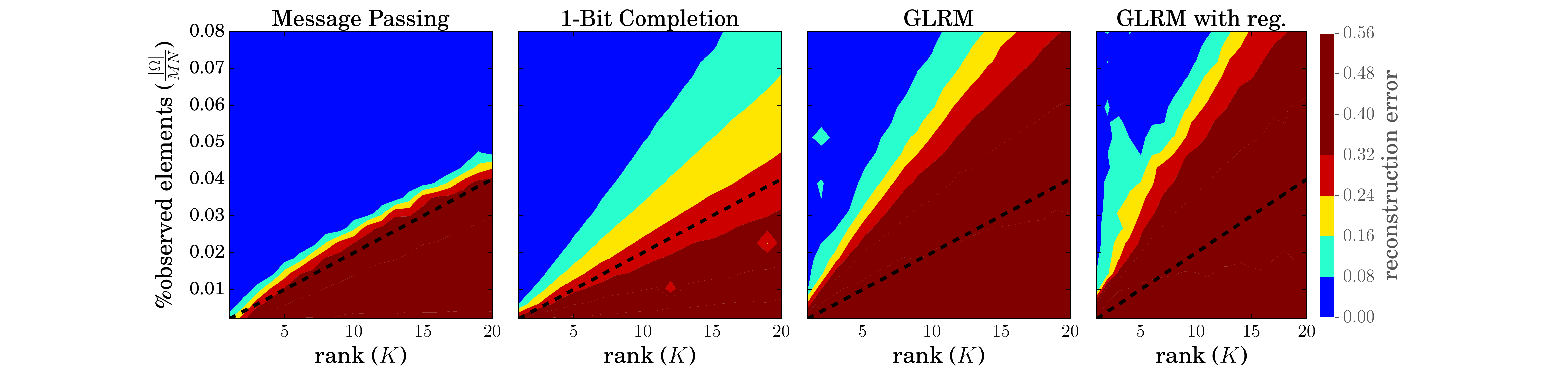}}
%\subcaptionbox*{}{\includegraphics[width=.3\textwidth]{figures/phase_glrm4x4.pdf}}
}
\vspace*{-.3in}%
  \caption{\it \small  
  The matrix completion error for Message Passing, 1-Bit matrix completion and GLRM (with and without regularization) 
as a function of matrix rank and portion of observed elements $|\Omega|$ for $\M = \N = 1000$. 
  The dashed black line indicates the tentative information bottleneck.}
  \label{fig:stuff}
\vspace*{-1em}
\end{figure*}

\section{Experiments}\label{sec:experiments}
We evaluated the performance of message passing 
on random matrices and real-world data. %\footnote{The Python implementation is available at \url{https://github.com/mravanba/BooleanFactorization}}
In all experiments,
message passing uses damping with $\lambda = .4$, 
$T = 200$ iterations and 
uniform priors $\PPX_{\m, \k}(1) = \PPY_{\k, \n}(1) = .5$.
This also means that if the channel is 
symmetric -- that is $\PPO(1 \mid 1) = \PPO(0 \mid 0) > .5$ -- the approximate MAP reconstruction $\widehat{\ZZ}$ does not depend on 
$\PPO$, and we could simply use $\PPO _{\m, \n}(1 \mid 1) = \PPO _{\m, \n}(1 \mid 1) = c $ for any $ c > .5$. The only remaining hyper-parameters 
are rank $K$ and maximum number of iterations $T$.

\begin{figure}
  \begin{center}
    \includegraphics[width=.4\textwidth]{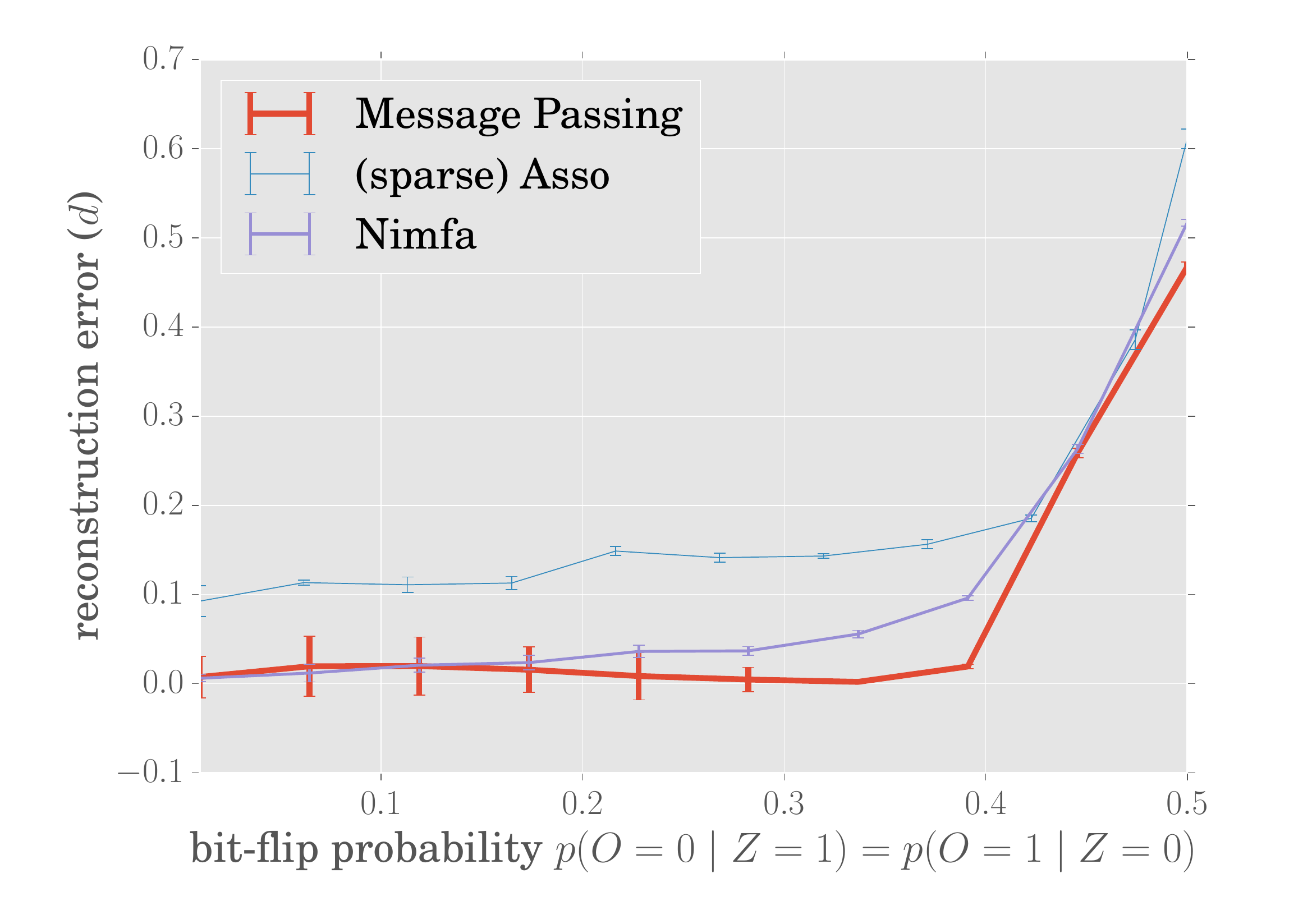}
  \end{center}
\vspace*{-2em}
  \caption{{\it \small Comparison of message passing and NIMFA for \textit{Boolean matrix factorization}}}\label{fig:nimfa}
\vspace*{-1em}
\end{figure}

\subsection{Random Matrices}
\noindent \textbf{Matrix Factorization.}
We compared our method against binary matrix factorization method of \cite{zhang2007binary}, which was implemented by NIMFA~\citep{ZitnikZ12} 
as well as (sparse) Asso of \citet{miettinen2006discrete}.
Here, all methods receive the correct $\K$ as input.

\Cref{fig:nimfa} compares the reconstruction error of different methods at different noise levels.
The results are for $1000 \times 1000$ random matrices of rank $\K = 5$ where $\XX$ and $\YY$ were uniformly sampled from binary matrices.
The results for different $K$ show a similar trend.%
\footnote{Both message passing and NIMFA use the same number of iterations $T = 200$. 
For NIMFA we use the 
default parameters of $\lambda_h = \lambda_w = 1.1$ and initialize the matrices using SVD. For Asso we report the result
for the best threshold hyper-parameter $\tau \in \{.10, .31, .52, .74, .95\}$.} 
The \textit{reconstruction error} is 
\begin{align}\label{eq:rec_error}
  \dd(\ZZ, \widehat{\ZZ}) \quad \defeq \quad \frac{1}{\M \N}\sum_{\m,\n} | \zZ_{\m,\n} - \widehat{\zZ}_{\m,\n}|.
\end{align}

The results suggests that message passing and NIMFA are competitive, with message passing performing better 
at higher noise levels.
The experiments were repeated $10$ times for each point. 
The small variance of message passing performance at low noise-levels
is due to the multiplicity of symmetric MAP solutions, and could be resolved by performing decimation,
albeit % RG
at a computational cost. 
We speculate that the symmetry breaking of higher noise levels help message passing choose a fixed point,
which results % ; resulting 
in lower variance. Typical running times for a single matrix in this setting are 2, 15 and 20 seconds for NIMFA, message passing and 
sparse Asso respectively.\footnote{Since sparse Asso is repeated 5 times for different hyper-parameters, its overall run-time is 100 seconds.}

Despite being densely connected, at lower levels of noise, BP often converges within the maximum number of iterations. 
The surprisingly good performance of BP, despite the large number of loops, is because most factors have a weak 
influence on many of their neighboring variables. This effectively limits the number of influential loops in the factor-graph; see Appendix C
for more.

\noindent \textbf{Matrix Completion.} The advantage of message passing to its competition is more evident in matrix ``completion'' problem,
where the complexity of BP grows with the number of observed elements, rather than the size of matrix $Z$.
%For $\K \ll \M,\N$, 
We can ``approximate'' a lower-bound on the \textit{number of observed entries} $|\Omega| = M N (1-\PPO(\nullmath))$
required for recovering $\ZZ$ by
\begin{align}
  \label{eq:bound}
|\Omega|\  >\  \K(\M + \N - \log(\K) + 1) + \OOO(\log(\K)).
\end{align}
To derive % arrive at [RG]
this approximation, we briefly sketch an information theoretic argument. 
Note that the total number of ways to define a Boolean matrix $\ZZ \in \{0,1\}^{\M \times \N}$ of rank $\K$ is 
$\frac{2^{\K (\M + \N)}}{\K!}$, where the nominator
is the number of different $\XX$ and $\YY$ matrices and $\K!$ is the irrelevant degree of freedom in choosing the permutation matrix $\UU$, such that
$\ZZ = (\XX \booltimes \UU) \booltimes (\UU^T \booltimes \YY)$.  
%\RG{Hmmm.. aren't there many (X, Y) pairs that produce the same Z ... Does this matter?}{}
The logarithm of this number, using Sterling's approximation, is the r.h.s.~of \Cref{eq:bound}, 
lower-bounding the number of bits required to recover $\ZZ$, in the absence of any noise. 
Note that this is assuming that any other degrees of freedom in producing $\ZZ$ grows
sub-exponentially with $K$ -- \ie is absorbed in the additive term $\OOO(\log(\K))$. 
%Although we cannot validate this assumption, we know the bound to be correct at least for $K=1$.
This approximation also resembles the $\OOO(\K \N \mathrm{polylog}(\N))$ sample complexity for various 
real-domain matrix completion tasks \citep[\eg][]{candes2010matrix,keshavan2010matrix}.

\Cref{fig:stuff} compares message passing against GLRM and 1-Bit matrix completion.
In all panels of \Cref{fig:stuff}, each point represents the average reconstruction error for random $1000 \times 1000$ Boolean matrices. 
For each choice of observation percentage $\frac{|\Omega|}{M N}$ and rank $\K$, the experiments were repeated $10$ times.%
\footnote{This means each figure summarizes $20\, \mbox{(rank)} \times 20\, \mbox{(number of observations)} \times 10\, \mbox{(repeats)} \,=\, 4000$
experiments. The exception is 1-Bit matrix completion, where due to its longer run-time the
number of repetition was limited to two. The results for 1-Bit completion are for best $\beta \in \{.1,1,10\}$.} 
The dashed black line is the information theoretic approximate lower-bound of \Cref{eq:bound}. 
This result suggests that message passing outperforms both of these methods and remains effective close to this bound.

%Therefore our experiment of \Cref{fig:stuff} also suggests that \textit{when the underlying factors are Boolean, using Boolean matrix completion works better than non-Boolean alternative in practice.}

\Cref{fig:stuff} also suggests that, when using message passing, the transition from recoverability to non-recoverability is sharp. 
Indeed the variance of the reconstruction error is always close to zero, but in a small neighborhood of the dashed black line.%
\footnote{The sparsity of $\ZZ$ is not apparent in  \Cref{fig:stuff}. 
Here, if we generate $\XX$ and $\YY$ uniformly at random, as 
 $\K$ grows, the matrix $\ZZ = \XX \booltimes \YY$ becomes all ones. To avoid this degeneracy, we choose $\PPX_{\m,\k}(\xX_{\m,\k})$ and $\PPY_{\k,\n}(\yY_{\k,\n})$ 
so as to enforce $\PP(\ZZ = 1) \approx \PP(\ZZ = 0)$. It is easy to check that $\PPX_{\m,\k}(1) = \PPY_{\k,\n}(1) = \sqrt{1 - \sqrt[\K]{.5}}$ produces this desirable outcome. Note that these probabilities are only used for random matrix ``generation'' and the message passing algorithm is
using uniform priors.}
%Typical run-time for message passing and GLRM is similar (\eg 5 seconds for $\K=10$ and $\frac{|\Omega|}{M N} = .02$), while 
%1-Bit completion is several times slower, taking 15 seconds in the same setting.

\subsection{Real-World Applications}
This section evaluates message passing on two real-world applications. 
While there is no reason to believe that the real-world matrices must necessarily decompose into low-rank Boolean factors, we see
that Boolean completion using message passing performs well in comparison with other methods that assume Real factors.

\begin{table}
  \caption{{\small \it Matrix completion performance for MovieLense dataset.} }
%1-Bit completion results are for best hyper-parameter $\beta \in \{.1, 1, 10\}$. All GLRM results are for best quadratic regularization
%hyper-parameter in $\{.01,.1,1,10\}$
%}.}
  \label{tab:lens}
  \centering
\resizebox{\columnwidth}{!}{%
  \begin{tabular}{c  l c c  c  c  c c  c c}
  && time (sec)&binary &\multicolumn{6}{c}{\textbf{observed} percentage of available ratings} \\
  && min-max   &input?&1\% & 5\% & 10\% & 20\% & 50\% & 95\%\\ \cline{2-10}\\\cline{2-10}
  \parbox[t]{2mm}{\multirow{3}{*}{\rotatebox[origin=c]{90}{\underline{1M-dataset}}}}  
                         & message passing & {2-43}&Y & \textbf{56\%} & 65\% & 67\% & 69\% & 71\% & 71\% \\\cline{3-10}
  &                       GLRM (ordinal hinge) & 2-141  & N & 48\% & 65\% & 68\%& 70\% & 71\%&  72\%\\\cline{3-10}
  &                       GLRM (logistic)      & 4-90 & Y & 46\% & 63\% & 63\%& 63\% & 63\% & 62\% \\\cline{2-10}\\\cline{2-10}

  \parbox[t]{2mm}{\multirow{4}{*}{\rotatebox[origin=c]{90}{\underline{100K-dataset}}}} & 
                          message passing & {0-2}& Y & \textbf{52}\% & \textbf{60}\% & \textbf{63}\% & 65\% & 67\% & 70\% \\\cline{3-10}
  &                       GLRM (ordinal hinge) & {0-2} & N & 48\% & 58\% & {63}\%& 67\% & 69\%&  70\%\\\cline{3-10}
  &                       GLRM (logistic) & {0-2} & Y & 45\% & 50\% & 62\%& 63\% & 62\%&  67\%\\\cline{3-10}
  &  1-bit completion & \textbf{30-500} &Y & 50\% & 53\% & 61\%& 65\% & {70}\%&  {72}\%\\\cline{2-10}
                         
  \end{tabular}
}
\end{table}

% \begin{table}
%   \caption{{\small \it Matrix completion performance for MovieLense-1M dataset}}
%   \label{tab:lens}
%   \centering
% \resizebox{\columnwidth}{!}{%
%   \begin{tabular}{c c | c  c  c c  c c}
%     &&\multicolumn{6}{c}{\textbf{observed} percentage of 1 million ratings} \\
%   && 1\% & 5\% & 10\% & 20\% & 50\% & 95\%\\ \cline{1-8}
%     {\multirow{2}{*}{K=2}} & Boolean factorization & \textbf{56\%} & 65\% & 67\% & 69\% & 71\% & 71\% \\\cline{2-8}
%                          & ordinal hinge loss & 48\% & 65\% & 68\%& 70\% & 72\%&  72\%\\\cline{1-8}
%     {\multirow{2}{*}{K=4}} & Boolean factorization & \textbf{55\%} & \textbf{61\%} & \textbf{63\%} & 65\% & 69\% & 70\% \\\cline{2-8}
%                          & ordinal hinge loss & 45\% & 50\% & 57\% & 62\% & 62\% & 71\% \\
%   \end{tabular}
% }
% \end{table}

\subsubsection{MovieLens Dataset}
We applied our message passing method to MovieLens-1M and MovieLens-100K dataset\footnote{\url{http://grouplens.org/datasets/movielens/}} as an application in \textit{collaborative filtering}.
The Movie-Lense-1M dataset contains 1 million ratings from 6000 users on 4000 movies (\ie $1/24$ of all the ratings are available). The ratings are ordinals 1-5. Here we say a user is ``interested'' in the movie iff her rating is above the global average of ratings. The task is to predict this single bit by observing a random subset of the available user$\times$movie rating matrix. For this, we use $\alpha \in (0,1)$ portion of the $10^6$ ratings to predict the one-bit interest level for the remaining ($1 - \alpha$ portion of the) data-points. Note that here $|\Omega| = \frac{\alpha \, M \, \N}{24}$.
The same procedure is applied to the smaller Movie-Lens-100K dataset. 
The reason for including this dataset was to compare message passing performance with 1-Bit matrix completion that
does not scale as well.

%For binary matrix completion, we threshold the training data to $\{0,1\}$ values from the start using the global average rating.
We report the results using GLRM with logistic and \textit{ordinal} hinge loss \citep{rennie2005loss} and quadratic regularization of the factors.
\footnote{The results reported for 1-Bit matrix completion are for best $\beta \in \{.1, 1, 10\}$ (see \Cref{eq:1bitcompletion}).
The results for GLRM are for the regularization parameter in $\{.01, .1,1,10\}$ with the best test error.}
%1-bit matrix completion as well as GLRM with logistic loss use thresholded ratings, 
Here, only GLRM with ordinal hinge loss uses actual ratings (non-binary) to predict the ordinal ratings which are then thresholded.
 
%For These ratings are then thresholded by the global average rating to predict the interest level.

\Cref{tab:lens} reports the run-time and test error of all methods for $K = 2$, using different 
$\alpha \in \{.01, .05,.1,.2,.5,.95\}$ portion of the available ratings. It is surprising that
only using one bit of information per rating, message passing and 1-bit completion are competitive with ordinal hinge loss that 
benefits from the full range of ordinal values. 
The results also suggest that when only few 
observations are available (\eg $\alpha  = .01$), message passing performs better than all other methods.
With larger number of binary observations, 1-bit completion performs slightly better than message passing, 
but it is orders of magnitude slower. Here, the variance in the range of reported times in \Cref{tab:lens}
is due to variance in the number of observed entries -- \ie $\alpha=.01$ often has the smallest run-time.%
%\footnote{\citet{davenport20141} report $73\%$ accuracy using 1-Bit matrix completion after observing $95\%$ of ratings in 100K-Movielens dataset.
%This slightly lower error is achieved by selecting the best $\beta$ in a larger set of 20 values.}

%Moreover, Boolean factorization seems to be less prone to over-fitting for $K = 20$. 
%In terms of running time, message passing, which was using a maximum of 100 iterations, was at least as fast as GLRM.%
%\footnote{On a 4GHz machine with 8 cores, message passing was 5-10 times faster.}

\begin{figure}
  \begin{center}
    \includegraphics[width=0.35\textwidth]{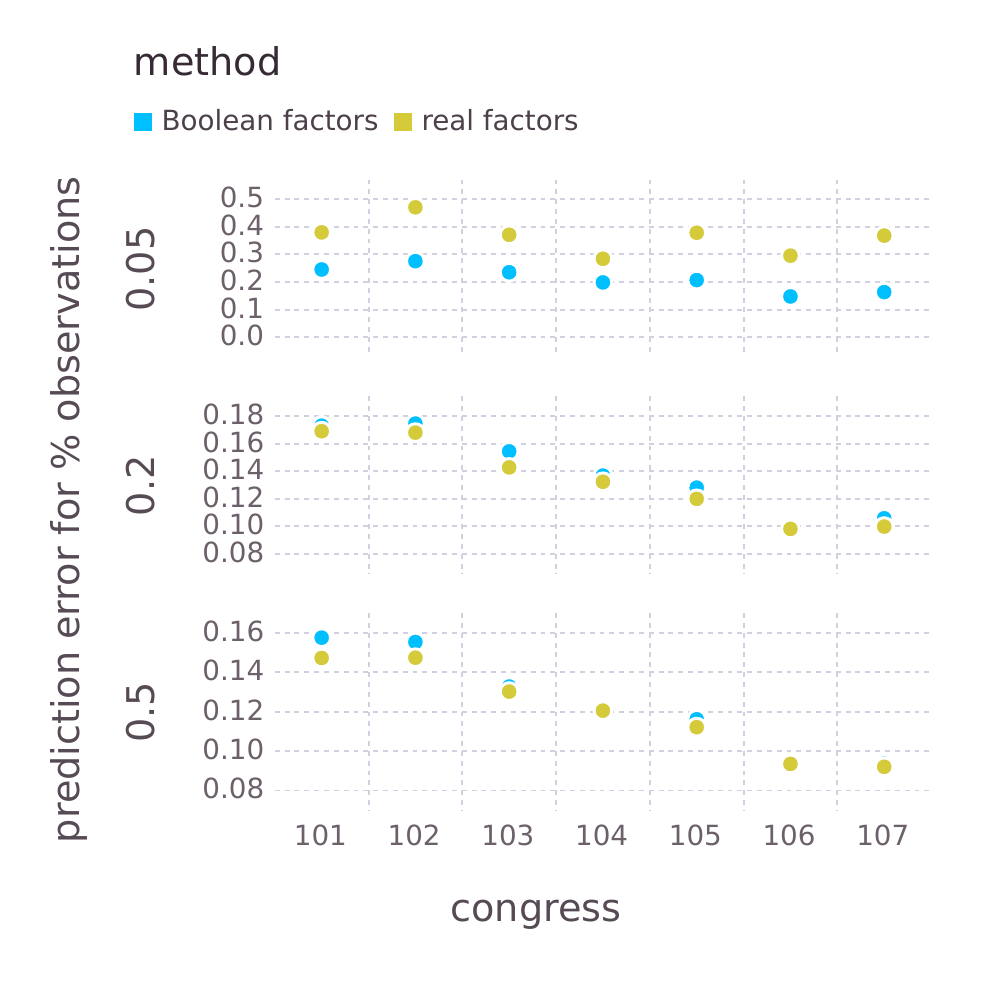}
  \end{center}
\vspace*{-2em}
  \caption{{\it \small The prediction error using Boolean matrix completion (by message passing) versus using GLRM with hinge loss for binary matrix completion using real factors. Each panel has a different observed percentage of entries $\frac{|\Omega|}{M N} \in \{.05,.2,.5\}$. Here the horizontal axis identifies senator$\times$issue  matrices and the y-axis is the average error
in prediction of the unobserved portion of the (yes/no) votes.}}\label{fig:senate}
\vspace*{-2em}
\end{figure}

% \citet{davenport20141} also report a good performance by only observing a single bit, instead of ordinal ratings for the MovieLens dataset, albeit on 
% a small version of this dataset (with 100,000 entries) and for the dense setting of  $\alpha = .95$.
% Similar to GLRM, they assume real factors, and according to their model, the 1-bit observation is obtained by thresholding a 
% function of the product matrix -- \ie $\OO = f(\ZZ) \overset{?}{\geq} c$ for some function $f$ and constant $c$.
% Using sigmoid function and $c = .5$, their method is closely related to GLRM with logistic loss (instead of the hinge loss, used in our experiments).

\subsubsection{Reconstructing Senate Voting Records}
We applied our noisy completion method to predict the (yes/no) senate votes during 1989-2003 by observing a randomly selected subset of votes.\footnote{The senate data was obtained from {\scriptsize\url{http://www.stat.columbia.edu/~jakulin/Politics/senate-data.zip}} 
prepared by \cite{jakulin2009analyzing}.}
 This dataset contains $7$ Boolean matrices (corresponding to voting sessions for $101^{st}-107^{th}$ congress), where a small portion  
of entries are missing. For example the first matrix is a $634 \times 103$ Boolean matrix recording the vote of 102 senators on $634$ topics plus the outcome of the vote (which we ignore).

\Cref{fig:senate} compares the prediction accuracy in terms of reconstruction error \Cref{eq:rec_error} of message passing and GLRM (with hinge loss or binary predictions) for the best choice of $\K \in \{1,\ldots,10\}$ on each of 7 matrices. \footnote{GLRM is using quadratic regularization  while message passing is using uniform priors.}
In each case we report the prediction accuracy on the unobserved entries, after observing $ \frac{|\Omega|}{M N} \in \{5\%, 20\%, 50\%\}$ of the votes. For sparse observations ($\frac{|\Omega|}{M N} = .05$), the message passing error is almost always half of the error when we use real factors. With larger number of observations, the methods are comparable, with GLRM performing slightly better.

\section*{Conclusion \& Future Work}
This paper introduced a simple message passing technique for approximate Boolean factorization and noisy matrix completion. While having a linear time complexity, 
%\RG{I thought it was NP-hard. Or do you mean EACH iteration ??}{}
 this procedure favorably compares with the state-of-the-art in Boolean matrix factorization and completion. 
In particular, for matrix completion with few entries, message passing significantly outperforms the existing methods that use real factors. 
This makes message passing a useful candidate for collaborative filtering in modern applications involving large datasets of sparse Boolean
observations.

Boolean matrix factorization with modular arithmetic, replaces the logical OR operation with exclusive-OR, only changing one of the factor types (\ie  type $\g$) in our graphical model. Therefore both min-sum and sum-product message passing can also be applied to this variation. The similarity of this type of Boolean factorization to LDPC codes,
suggests that
one may be able to use noisy matrix completion as an efficient method of communication over a noisy channel, where the data is preprocessed to have low-rank matrix form and a few of its entries are then transmitted through the noisy channel. This is particularly interesting, as both the code and its parity checks are transmitted as a part of the same matrix. We leave this promising direction to future work.

\bibliographystyle{plainnat}
\bibliography{refs.bib}

\clearpage
\appendix
\section{Detailed Derivation of Simplified BP Messages}

The sum of the factors in the factor-graph of \Cref{fig:fg} is
\begin{align}
  &\sum_{\m,\k} \h_{\m,\k}(\xX_{\m,\k}) + \sum_{\n,\k} \h_{\n,\k}(\yY_{\n,\k})  + \notag \\ 
  & \sum_{\m,\n,\k} \f_{\m,\n,\k}(\xX_{\m,\n,\k}, \yY_{\m,\n,\k}, \wW_{\m,\n,\k}) + \notag \\ 
&\sum_{\m,\n} \g_{\m,\n}(\{\wW_{\m,\n,\k} \}_\k) \label{eq:sum_factors} \\
=&\sum_{\m,\n} \log(\PPX(\xX_{\m,\k})) + \sum_{\n,\k} \log(\PPY(\yY_{\k,\n})) + \notag \\
&\sum_{\m,\n,\k} \ident(\wW_{\m,\n,\k} = \xX_{m,\k} \wedge \yY_{\k,\n}) + \notag \\
&\sum_{\m,\n} \log\big(\PPO_{\m,\n}(\oO_{\m,\n} \mid \bigvee_{\k} \wW_{\m,\n,\k})\big) \label{eq:replace}\\
=& \sum_{\m,\n} \log(\PPX(\xX_{\m,\k})) + \sum_{\n,\k} \log(\PPY(\yY_{\k,\n})) + \notag \\
&\sum_{\m,\n} \log\big(\PPO_{\m,\n}(\oO_{\m,\n} \mid \bigvee_{\k} \xX_{m,\k} \wedge \yY_{\k,\n})\big) \label{eq:absorb}\\
= &\log(\PP(\XX, \YY \mid \OO)) \label{eq:logposterior}
\end{align}
where in \Cref{eq:replace} we replaced each factor with its definition. \Cref{eq:absorb}  combines the two last terms of \Cref{eq:replace},
which is equivalent to marginalizing out $\WW$. The final result of \Cref{eq:logposterior} is the log-posterior of \Cref{eq:map}.

Since the original MAP inference problem of $\arg_{\XX, \YY} \max \; \PP(\XX, \YY \mid \OO)$ is equivalent to $\arg_{\XX, \YY} \max \; \log(\PP(\XX, \YY \mid \OO))$,
our objective is to perform max-sum inference over this factor-graph, finding an assignment that maximizes the summation of \Cref{eq:sum_factors}

We perform this max-sum inference using Belief Propagation (BP).  Applied to a factor-graph, BP involves message exchange between neighboring variable and factor nodes.
Two most well-known variations of BP are sum-product BP for marginalization and max-product or max-sum BP for MAP inference.
Here, we provide some details on algebraic manipulations that lead to the simplified form of max-sum BP message updates of \Cref{eq:map_update}.
\Cref{sec:var2fac} obtains the updates \Cref{eq:xout} and \Cref{eq:yout} in our algorithm and
\Cref{sec:fac2var} reviews the remaining message updates of \Cref{eq:map_update}

\subsection{Variable-to-Factor Messages}\label{sec:var2fac}

Consider the binary variable $\xX_{\m, \k} \in \{0,1\}$ in the graphical model of \Cref{fig:fg}.
Let $\msgt{\xX_{\m, \k}}{\f_{\m,\n,\k}}(\xX_{\m, \k}): \{0,1\} \to \Re$ be the \textit{message} from variable $\xX_{\m, \k}$ to 
the factor $\f_{\m,\n,\k}$ in this factor-graph. Note that this message contains two assignments for $\xX_{\m, \k} = 0$ and $\xX_{\m, \k} = 1$.
As we show here, in our simplified updates this message is represented by $\aAh_{\m,\n,\k}$.
In the max-sum BP, the outgoing message from any variable to a neighboring factor is the sum of all incoming messages,
except for the message from the receiving factor -- \ie
\begin{align}
 & \msgt{\xX_{\m, \k}}{\f_{\m,\n,\k}}(\xX_{\m, \k})\ttt{t+1} = \msgt{\h_{\m,\k}}{\xX_{\m, \k}}(\xX_{\m, \k})\ttt{t}  \notag \\
 & + \, \sum_{\n' \neq \n} \msgt{\f_{\m,\n',\k}}{\xX_{\m, \k}}(\xX_{\m, \k})\ttt{t} \; + \; c \label{eq:bp_var}
\end{align}

What matters in BP messages is the difference between the message $\msgt{\xX_{\m, \k}}{\f_{\m,\n,\k}}(\xX_{\m, \k})$ assignment  
for $\xX_{\m, \k} = 1$ and $\xX_{\m, \k} = 0$ (note the constant $c$ in \Cref{eq:bp_var}). Therefore we can use a singleton message value that capture this difference instead of using a message over the
binary domain -- \ie
\begin{align}\label{eq:bias}
\aAh_{\m,\n,\k} = \msgt{\xX_{\m, \k}}{\f_{\m,\n,\k}}(1) - \msgt{\xX_{\m, \k}}{\f_{\m,\n,\k}}(0)
\end{align}

This is equivalent to assuming that the messages are normalized so that $\msgt{\xX_{\m, \k}}{\f_{\m,\n,\k}}(0) = 0$. We will extensively use this normalization assumption in the following. 
By substituting \Cref{eq:bp_var} in \Cref{eq:bias} we get the simplified update of \Cref{eq:xout}
\begin{align*}
  \aAh_{\m,\n,\k}\ttt{t+1} &= 
    \bigg ( \msgt{\h_{\m,\k}}{\xX_{\m, \k}}(1)\ttt{t} \; + \; \sum_{\n' \neq \n} \msgt{\f_{\m,\n',\k}}{\xX_{\m, \k}}(1)\ttt{t}(1) \bigg ) \\
&-  \bigg ( \msgt{\h_{\m,\k}}{\xX_{\m, \k}}(0)\ttt{t} \; +  \sum_{\n' \neq \n} \msgt{\f_{\m,\n',\k}}{\xX_{\m, \k}}(0)\ttt{t}
\bigg ) \\
& = \bigg (\msgt{\h_{\m,\k}}{\xX_{\m, \k}}(1)\ttt{t} - \msgt{\h_{\m,\k}}{\xX_{\m, \k}}(0)\ttt{t}\bigg) \\
&+ 
 \sum_{n' \neq n} \bigg ( \msgt{\f_{\m,\n',\k}}{\xX_{\m, \k}(1)\ttt{t}} - \msgt{\f_{\m,\n',\k}}{\xX_{\m, \k}}(0)\ttt{t} \bigg ) \\
& = \log \bigg (\frac{\PPX_{\m,\k}(1)}{\PPX_{\m,\k}(0)} \bigg) + \sum_{\n' \neq \n} \; \aA_{\m,\n',\k}\ttt{t}
\end{align*}
 and we used the
fact that %$\h_{\m,\k}$ and $\f_{\m,\n,\k}$ are
\begin{align*}
  &\aA_{\m,\n',\k} =  {\msgt{\f_{\m,\n',\k}}{\xX_{\m, \k}}(1)\ttt{t}} \; - \; {\msgt{\f_{\m,\n',\k}}{\xX_{\m, \k}}(0)\ttt{t}} \\
  &\log \bigg (\frac{\PPX_{\m,\k}(1)}{\PPX_{\m,\k}(0)} \bigg) = \h_{\m,\k}(1) - \h_{\m,\k}(0).
\end{align*}

The messages $\bBh_{\m,\n,\k}$ from the variables $\yY_{\n,\k}$ to $\f_{\m,\n,\k}$ is obtain similarly. 
The only remaining variable-to-factor messages in the factor-graph of \Cref{fig:fg} are from auxiliary variables $\wW_{\m,n,k}$ to neighboring factors.
However, since each variable $\wW_{\m,\n,k}$ has exactly two neighboring factors, the message from $\wW_{\m,\n,\k}$ to any
of these factors is simply the incoming message from the other factor -- that is
\begin{align}
  &\msgt{\wW_{\m,\n,\k}}{\g_{\m,\n}}(\wW_{\m,\n,\k}) = \msgt{\f_{\m,\n,\k}}{\wW_{\m,\n,\k}}(\wW_{\m,\n,\k}) \notag \\
&  \msgt{\g_{\m,\n}}{\wW_{\m,\n,\k}}(\wW_{\m,\n,\k}) = \msgt{\wW_{\m,\n,\k}}{\f_{\m,\n,\k}}(\wW_{\m,\n,\k}) \label{eq:local_out}
\end{align}

\subsection{Factor-to-Variable Messages}\label{sec:fac2var}
The factor-graph of \Cref{fig:fg} has three types of factors. We obtain the simplified messages from each of these factors to their neighboring variables in the following sections.

\subsubsection{Local Factors} 
The local factors are $\{\h_{\m,k}\}_{\m,\k}$ and $\{\h_{n,k}\}_{\n,\k}$, each of which
is only connected to a single variable. The unnormalized message, leaving these factors is identical to the factor itself. 
We already used the normalized messages from these local factors
to neighboring variables in \Cref{eq:local_out} -- \ie $\h_{\m,k}(1) - \h_{\m,k}(0)$ and $\h_{\n,k}(1) - \h_{\n,k}(0)$, respectively.

\subsubsection{Constraint Factors} 
The constraint factors $\{\f_{\m,\n,\k}\}_{\m,\n,\k}$ ensure $\forall_{\m,\n,\k} \wW_{\m,\n,\k} = \xX_{\m,\k} \wedge \yY_{\n,\k}$.
Each of these factors has three neighboring variables.
In max-sum BP the message from a factor to a neighboring variable is given by the sum of that factor and incoming messages from its neighboring variables, except for the receiving variable, max-marginalized
over the domain of the receiving variable. Here we first calculate the messages from a constraint factor to $\xX_{\m,\k}$ (or equivalently $\yY_{\n,\k}$) variables in \textbf{(1)}. In \textbf{(2)} we derive the simplified messages to the auxiliary variable $\wW_{\m,\n,\k}$.

\noindent \textbf{(1)} according to max-sum BP equations the message from the factor $\f_{\m,\n,\k}$ to variable ${\xX_{\m,\k}}$ is 
\begin{align*}
&\msgt{\f_{\m,\n,\k}}{\xX_{\m,\k}}(\xX_{\m,\k})\ttt{t+1} = \\
&\max_{\wW_{\m,\n,\k}, \yY_{\n,\k}} \bigg ( \f_{\m,\n,\k}(\xX_{\m,\k}, \wW_{\m,\n,\k},\yY_{\n,\k})  \\ &+ \, \msgt{\yY_{\n,\k} \f_{\m,\n,\k}}(\yY_{\n,\k})\ttt{t}  
+\, \msgt{\wW_{\m, \n,\k}}{\f_{\m,\n,\k}}(\wW_{\m, \n,\k})\ttt{t} \bigg )
\end{align*}

For notational simplicity we temporarily use the shortened version of the above
\begin{align} \label{eq:max_aux_x}
  \msg_1'(\xX) = \max_{\wW, \yY} \f(\xX, \wW, \yY) \;
  + \msg_2(\yY)  +  \msg_3(\wW) 
\end{align}
where 
\begin{align*}
&  \msg_1(\xX) = \msgt{\xX_{\m,\k}}{\f_{\m,\n,\k}}(\xX_{\m,\k}) \\ 
& \msg'_1(\xX) = \msgt{\f_{\m,\n,\k}}{\xX_{\m,\k}}(\xX_{\m,\k}) \\
&  \msg_2(\yY) = \msgt{\yY_{\n,\k}}{\f_{\m,\n,\k}}(\yY_{\n,\k}) \\ 
& \msg'_2(\yY) = \msgt{\f_{\m,\n,\k}}{\yY_{\n,\k}}(\yY_{\n,\k}) \\
&  \msg_3(\wW) = \msgt{\wW_{\m,\n,\k}}{\f_{\m,\n,\k}}(\wW_{\m,\n,\k}) \\
& \msg'_3(\wW) = \msgt{\f_{\m,\n,\k}}{\wW_{\m,\n,\k}}(\wW_{\m,\n,\k}),
\end{align*}
that is we use $\msg(.)$ to denote the incoming messages to the factor and $\msg'(.)$ to identify the outgoing message.

If the constraint $\f(\xX, \yY, \wW) = \ident(\wW = \xX \wedge \yY) $ is not satisfied by an assignment to $\xX, \yY$ and $\wW$, it evaluates to $-\infty$, and therefore it does not have any effect on the outgoing message
due to the $\max$ operation. Therefore we should consider the $\max_{\wW, \yY}$ only over the assignments that satisfy $\f(.)$.

Here, $\xX$ can have two assignments; for $\xX = 1$, if $\yY = 1$, then $\wW = 1$ is enforced by $\f(.)$, and if
$\yY = 0$ then $\wW = 0$. Therefore \Cref{eq:max_aux_x} for $\xX = 1$ becomes
\begin{align}\label{eq:f2x1}
  \msg_1'(1) = \max ( \msg_2(1) + \msg_3(1), \msg_2(0) + \msg_3(0) )
\end{align}

For $\xX = 0$, we have $\wW = 0$, regardless of $\yY$ and the update of \Cref{eq:max_aux_x} reduces to
\begin{align} \label{eq:f2x2}
  \msg'_1(0) \; & = \; \max ( \msg_2(1) + \msg_3(0), \msg_2(0) + \msg_3(0) \} \\
&= \msg_3(0) + \max \{ \msg2(0), \msg2(1) \} \notag
\end{align}

Assuming the incoming messages are normalized such that $\msg_3(0) = \msg_2(0) = 0$ and denoting 
$$\bBh_{\m,\n,\k} = \msgt{\yY_{\n,\k}}{\f_{\m,\n,\k}}(1) - \msgt{\yY_{\n,\k}}{\f_{\m,\n,\k}}(0) =  \msg_2(1)$$ 
and 
$$\cCh_{\m,\n,\k} = \msgt{\wW_{\m,\n,\k}}{\f_{\m,\n,\k}}(1) - \msgt{\wW_{\m,\n,\k}}{\f_{\m,\n,\k}}(0) =  \msg_3(1)$$
the difference of \Cref{eq:f2x1} and \Cref{eq:f2x2} gives the normalized outgoing message of \Cref{eq:xin}
\begin{align}
  \aA_{\m,\n,\k} = &\msg'_1(1) - \msg'_1(0) = \max ( \cCh_{\m,\n,\k} + \bBh_{\m,\n,\k}, 0  ) \notag \\
- &\max (0,\bBh_{\m,\n,\k}) \label{eq:w2x_final}
\end{align}
The message of \Cref{eq:yin} from the constraint $\f_{\m,\n,\k}$ to $\yY_{\n, \k}$ is obtained in exactly the same way.\\

\noindent \textbf{(2)} The max-sum BP message from the constraint factor $\f_{\m,\n,\k}$ to the auxiliary variable $\wW_{\m,\n,\k}$ is
\begin{align*}
&\msgt{\f_{\m,\n,\k}}{\wW_{\m,\n,\k}}(\wW_{\m,\n,\k})\ttt{t+1} = \\
&\max_{\xX_{\m,\k}, \yY_{\n,\k}} \big ( \f_{\m,\n,\k}(\xX_{\m,\k}, \wW_{\m,\n,\k}, \yY_{\n,\k})  + \\
&  \msgt{\yY_{\n,\k}}{\f_{\m,\n,\k}}(\yY_{\n,\k})\ttt{t}  +  \msgt{\xX_{\m,\k}}{ \f_{\m,\n,\k}}(\wW_{\m,\n,\k})\ttt{t} \big )
\end{align*}

Here, again we use the short notation 
\begin{align} \label{eq:max_aux_w}
  \msg'_3(\wW) = \max_{\xX, \yY} \f(\xX, \wW, \yY) \;+
  \msg_1(\xX) + \msg_2(\yY)
\end{align}
and consider the outgoing message $\msg'(\wW)$ for $\wW = 1$ and $\wW = 0$.
If $\wW = 1$, we know that $\xX = \yY = 1$. This is because otherwise the factor $\f$ evaluates to $-\infty$. This simplifies \Cref{eq:max_aux_w} to
\begin{align*}
\msg_3'(1) = \msg_1(1) + \msg_2(1)
\end{align*}

For $\wW = 0$, either $\xX = 0$, or $\yY = 0$ or both. This means
\begin{align*}
  \msg'_3(0) = \max (&\msg_1(0) + \msg_2(1), \msg+1(1) + \msg_2(0), \\
&\msg_1(0) + \msg_2(0))
\end{align*}

Assuming the incoming messages were normalized, such that $\msg_2(0) = \msg_1(0) = 0$, the 
normalized outgoing message $\cC_{\m,\n,\k} = \msg_3(1) - \msg_3(0)$ simplifies to 
\begin{align}
  \cC_{\m,\n,\k} &= \msg_1(1) + \msg_2(1) - \max ( 0, \msg_1(1), \msg_2(1) ) \notag \\
  & = \min ( \msg_1(1) + \msg_2(1) , \msg_1(1), \msg_2(1) ) \notag \\
  & = \min ( \aAh_{\m,\n,\k} + \bBh_{\m,\n,\k},  \aAh_{\m,\n,\k}, \bBh_{\m,\n,\k} ) \notag
\end{align}

\subsection{Likelihood Factors}
At this point we have derived all simplified message updates of \Cref{eq:map_update}, except for the message $\cCh_{\m,\n,\k}$ from
factors $\g_{m,n}$ to the auxiliary variables $\wW_{\m,\n,\k}$ (\Cref{eq:zin}). These factors encode the likelihood term in the
factor-graph.% and deriving the message update for them is slightly more involved. 

% Recall that
% \begin{align}
%   \g_{\m,\n}(\{\wW_{\m,\n,\k}\}_\k) \; = \; 
% \log\big(\PPO_{\m,\n}(\oO_{\m,\n} \mid \bigvee_{\k} \wW_{\m,\n,\k})\big)% \; = \; \oR_{\m,\n}\; (-1)^{1 + (\bigvee_{\k} \wW_{\m,\n,\k})}
% \end{align}

The naive form of max-sum BP for the messages leaving this factor to each of $\K$ neighboring variables $\{\wW_{\m,\n,\k}\}_{1 \leq \k \leq \K}$ is
\begin{align} \label{eq:bp_g}
&\msgt{\g_{\m,\n}}{\wW_{\m,\n,\k}}(\wW_{\m,\n,\k})\ttt{t+1} = \\ 
&\max_{\{\wW_{\m,\n,\ell}\}_{\ell \neq \k \}}} \bigg ( \g_{\m,\n}(\{\wW_{\m,\n,{\ell'}}\}_{\ell'})  + \notag \\
&\sum_{\k' \neq \k} \msgt{\wW_{\m,\n,\k'}}{\g_{\m,\n}}(\wW_{\m,\n,\k'})\ttt{t} \bigg )\notag
\end{align}

However, since $\g(.)$ is a high-order factor (\ie depends on many variables), this naive update
has an exponential cost in $\K$. Fortunately, by exploiting the special form of $\g(.)$, we can reduce this cost to linear in $\K$.

In evaluating $\g(\{\wW_{\m,\n,\k}\}_\k)$ two scenarios are conceivable:
\begin{enumerate}
\item at least one of $\wW_{\m,\n,1},\ldots,\wW_{\m,\n,\K}$ is non-zero -- that is $\bigvee_{\k} \wW_{\m,\n,\k} = 1$ and $\g(\wW_{\m,\n,\k})$ evaluates to $\PPO_{\m,\n}(\oO_{\m,\n} \mid 1)$. 
\item $\bigvee_{\k} \wW_{\m,\n,\k} = 0$ and $\g(\wW_{\m,\n,\k})$ evaluates to $\PPO_{\m,\n}(\oO_{\m,\n} \mid 0)$.
\end{enumerate}
We can divide the maximization of \Cref{eq:bp_g} into two separate maximization operations over sets of assignments depending on the conditioning above and select the maximum of the two.

For simplicity, let $\msg_1(\wW_1),\ldots,\msg_\K(\wW_\K)$ denote 
$\msgt{\wW_{\m,\n,1}}{\g_{\m,\n}}(\wW_{\m,\n,1})\ttt{t},\ldots,\msgt{\wW_{\m,\n,\K}}{\g_{\m,\n}}(\wW_{\m,\n,\K})\ttt{t}$ 
respectively. W.L.O.G., let us assume the objective is to calculate the outgoing message to the first variable $\msg'_1(\wW_1) = \msgt{\g_{\m,\n}}{\wW_{\m,\n,1}}(\wW_{\m,\n,1})\ttt{t+1}$.
Let us rewrite \Cref{eq:bp_g} using this notation: 
\begin{align*}
\msg'_1(\wW_1) = \max_{\wW_{2} \ldots \wW_{\K}}  \big ( \g_{\m,\n}(\{\wW_\k\})  + 
  \sum_{\k' > 1} \msg_{\k'}(\wW_{\k'}) \big )
\end{align*}

For $\wW_1 = 1$, regardless of assignments to $\wW_2,\ldots,\wW_\K$, we have $\bigvee_{\k} \wW_{\m,\n,\k} = 1$ and therefore the maximization
above simplifies to
\begin{align*}
  \msg'_1(1) &= \max_{\wW_{2} \ldots \wW_{\K}}  \big ( \log(\PPO_{\m,\n}(\oO_{\m,\n} \mid 1))
  \sum_{\k' > 1} \msg_{\k'}(\wW_{\k'}) \big ) \notag \\
&= \log(\PPO_{\m,\n}(\oO_{\m,\n} \mid 1)) + \sum_{\k' > 1} \max (\msg_{\k'}(0), \msg_{\k'}(1)).
\end{align*}

For $\wW_1 = 0$, if $\forall_{\k' > 1} \wW_{\k'} = 0$ then $\g(\{\wW_\k\})$ evaluates to $\log(\PPO_{\m,\n}(\oO_{\m,\n} \mid 0)$, and otherwise
it evaluates to $\log(\PPO_{\m,\n}(\oO_{\m,\n} \mid 1)$. We need to choose the maximum over these two  cases. Note that in the second case we have to ensure at least one of the remaining variables is non-zero -- \ie $\exists_{\k' > 1} \wW_{\k'} = 1$. In the following update to enforce this constraint we use 
\begin{align}
  \label{eq:kstar}
\k^* = \arg_{\k' > 1} \max \msg_{\k'}(1) - \msg_{\k'}(0)  
\end{align}
to get
\begin{align*}
\msg'_1(0) = \max \bigg ( 
&\log(\PPO_{\m,\n}(\oO_{\m,\n} \mid 0) + \sum_{\k' > 1} \msg_{\k'}(0) \;,\; \\
&\log(\PPO_{\m,\n}(\oO_{\m,\n} \mid 1) + \msg_{\k^*} + \\
&\sum_{\k' > 1, \k' \neq \k^*} \max (\msg_{\k'}(0), \msg_{\k'}(1) )
\bigg )
\end{align*}
where, choosing $\wW_{\k^*} = 1$ maximizes the second case (where at least one $\wW_{\k'}$ for $\k' > 1$ is non-zero).

As before, let us assume that the incoming messages are normalized such that $\forall _{\k'} \msg_{\k'}(0) = 0$, and therefore $\cC_{\m,\n,\k'} = \msg_{\k'}(1)$. The normalized outgoing message is 
\begin{align*}
&\cCh_{\m,\n,1} = \msg'_1(1) - \msg'(0) = \log(\PPO_{\m,\n}(\oO_{\m,\n} \mid 1)) \\
&+ \sum_{\k' > 1} \max (0, \msg_{\k'}(1)) - \\
& \max \bigg ( \log(\PPO_{\m,\n}(\oO_{\m,\n} \mid 0) , \log(\PPO_{\m,\n}(\oO_{\m,\n} \mid 1) + \msg_{\k^*} \\
&+ \sum_{\k' > 1, \k' \neq \k^*}  \max (0, \msg_{\k'}(1) )
\bigg )\\
&= \min \bigg (\log(\PPO_{\m,\n}(\oO_{\m,\n} \mid 1)) - \log(\PPO_{\m,\n}(\oO_{\m,\n} \mid 0) \\
&+ \sum_{\k' > 1} \max (0, \msg_{\k'}(1)),
\max(-\msg_{\k^*}(1),0) \bigg )\\
& = \min \bigg ( \sum_{\k' > 1} \max(0, \cC_{\m,\n,\k'}\ttt{t}) \\
&+ \log \bigg(\frac{\PPO_{\m,\n}(\oO_{\m,\n} \mid 1)}{\PPO_{\m,\n}(\oO_{\m,\n} \mid 0)} \bigg ) , \max(0, -\max_{\k > 1}\cC_{\m,\n,\k'}\ttt{t})\bigg)
\end{align*}
where in the last step we used the definition of factor $\g$ and \Cref{eq:kstar} that defines $\msg_{\k^*}(1)$. This produces the simplified form of BP messages for the update \Cref{eq:zin} in our algorithm.
%\end{document}

\section{Marginal-MAP}\label{sec:marginal}
While the message passing for MAP inference approximates the ``jointly'' optimal assignment to $\XX$ and $\YY$ in the Bayesian setting, the marginals 
$\PP(\xX_{\m,\k} \mid \OO)$ and $\PP(\xX_{\k,\n} \mid \OO)$ are concerned with optimal assignments to ``individual'' $\xX_{\m, \k}$ and $\yY_{\k, \n}$ for each $\m,\n$ and $\k$. 
Here again, message passing can approximate the log-ratio of these marginals. 

We use the function $\lpx(a)= \log(1 + \exp(a))$ and its inverse $\lpx^{-1}(b) = \log(\exp(b) - 1)$ in the following updates for marginalization.
%\begin{subequations}
%\label{eq:marg_update}
\begin{align*}
&\aA_{\m,\n,\k}\ttt{t+1}  :=  {\cCh_{\m,\n,\k}\ttt{t} + \bBh_{\m,\n,\k}\ttt{t}} \, - \\
& \, \log \big(1 + \exp(\bBh_{\m,\n,\k}\ttt{t}) + \exp(\aAh_{\m,\n,\k}\ttt{t}) \big) \\
&\bB_{\m,\n,\k}\ttt{t+1}  :=  {\cCh_{\m,\n,\k}\ttt{t} + \aAh_{\m,\n,\k}\ttt{t}} \, - \\
& \, \log \big(1 + \exp(\aAh_{\m,\n,\k}\ttt{t} + \exp(\bBh_{\m,\n,\k}\ttt{t}) \big) \\
&\aAh_{\m,\n,\k}\ttt{t+1} :=  \log \bigg (\frac{\PPX_{\m,\k}(1)}{\PPX_{\m,\k}(0)} \bigg) + \sum_{\n' \neq \n} \; \aA_{\m,\n',\k}\ttt{t}  \\
&\bBh_{\m,\n,\k}\ttt{t+1} :=   \log \bigg(\frac{\PPY_{\n,\k}(1)}{\PPY_{\n,\k}(0)} \bigg ) + \sum_{\m' \neq \m} \; \bB_{\m',\n,\k}\ttt{t}  \\
&\cC_{\m,\n,\k}\ttt{t+1}  :=  \aAh_{\m,\n,\k}\ttt{t} + \bBh_{\m,\n,\k}\ttt{t} \\
&\cCh_{\m,\n,\k}\ttt{t+1} :=  \sum_{\k' \neq \k} \lpx(\cC_{\m,\n,\k'}\ttt{t}) \, + \, \log \bigg(\frac{\PPO_{\m,\n}(\oO_{\m,\n} \mid 1)}{\PPO_{\m,\n}(\oO_{\m,\n} \mid 0)} \bigg ) \\
&- \, \lpx \bigg( \phi^{-1} \big (\sum_{\k' \neq \k} \lpx(\cC_{\m,\n,\k'}\ttt{t}) \big) \\
&+ \, \log\bigg(\frac{\PPO_{\m,\n}(\oO_{\m,\n} \mid 1)}{\PPO_{\m,\n}(\oO_{\m,\n} \mid 0)} \bigg) \bigg)
\end{align*}
%\end{subequations}

Here, again using \Cref{eq:marg}, we can recover $\XX$ and $\YY$ from the marginals. 
However, due to the symmetry of the set of solutions,
one needs to perform \textit{decimation} to obtain an assignment to $\XX$ and $\YY$. Decimation 
 is the iterative process of running message passing then fixing the most biased variable --
\eg an $\xX_{\m,\k}\ \in\ \arg_{\m,\k}\max |\margX_{\m,\k}|$ -- after each convergence.
While a simple  randomized initialization of messages is often enough to break the symmetry of the solutions in max-sum inference,
in the sum-product case one has to repeatedly fix a new subset of most biased variables.

\begin{figure*}
  \centering
\hbox{
  \subcaptionbox*{}{\rotatebox[origin=t]{90}{\tiny \hspace{.08in}  observed  \hspace{.2in}   original  }}
  \subcaptionbox*{}{\includegraphics[width=.2\textwidth,angle=-90]{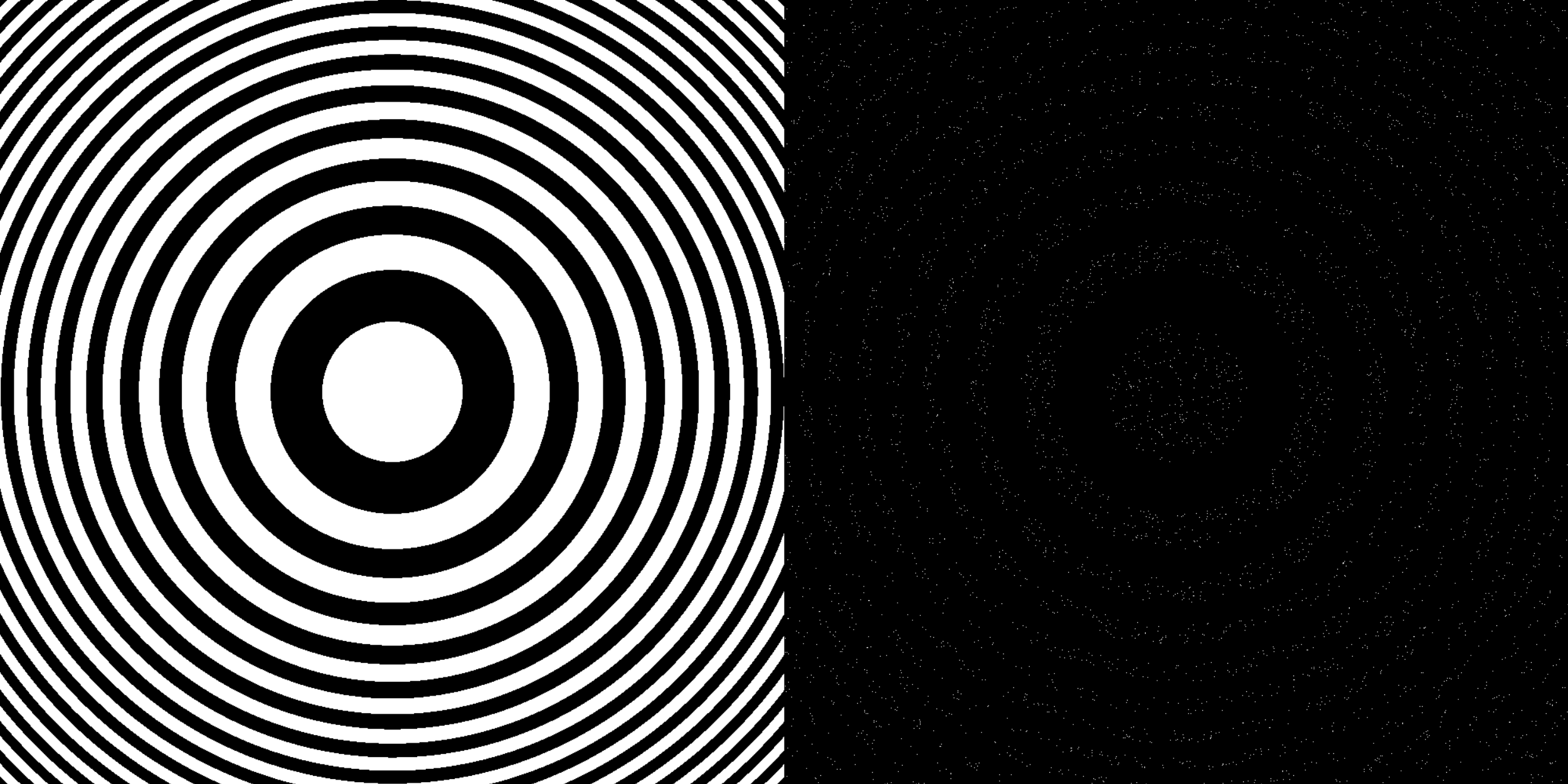}}
\hspace{.05in}
  \subcaptionbox*{}{\rotatebox[origin=t]{90}{\tiny \hspace{.05in} GLRM  \hspace{.2in}   message passing}}
  \subcaptionbox{\small $\rho \in \{.01, \,.02,\,.05\}$, $K = 10$}
{\includegraphics[width=.2\textwidth,angle=-90]{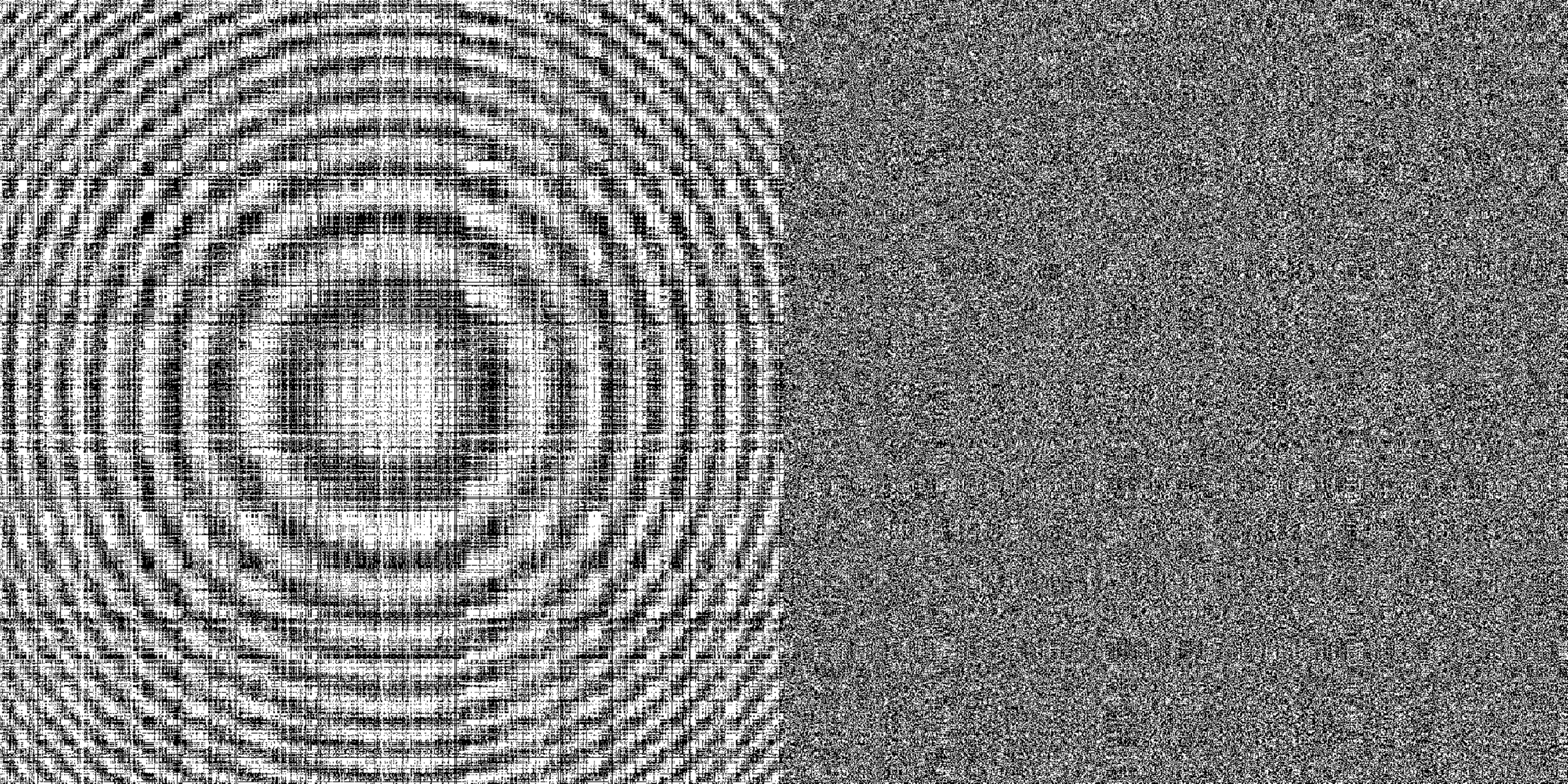}
  %\subcaptionbox{\small $\rho = .02$}
\includegraphics[width=.2\textwidth,angle=-90]{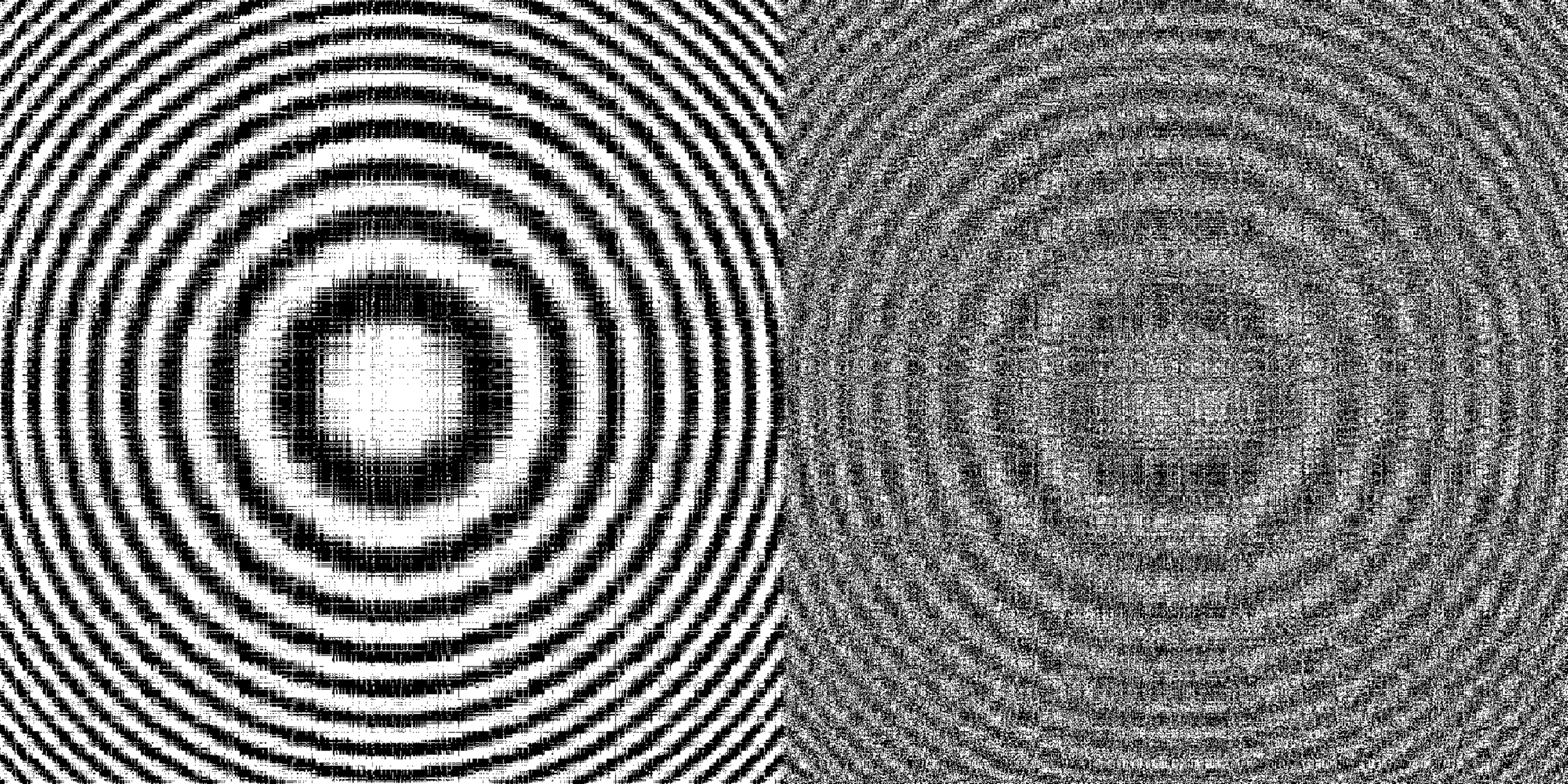}
%  \subcaptionbox{\small $\rho = .05$}
\includegraphics[width=.2\textwidth,angle=-90]{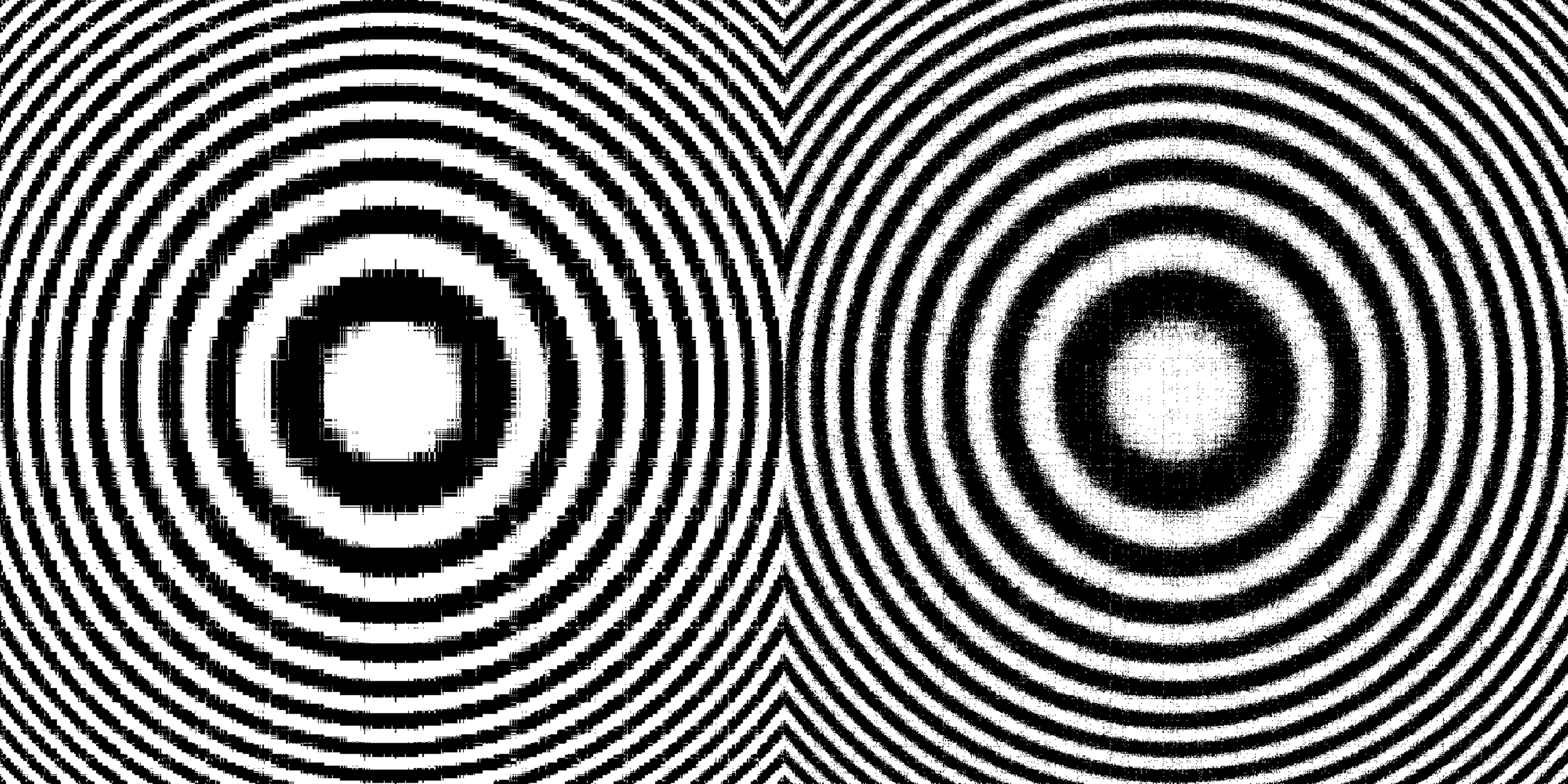}}
\hspace{.05in}
  %\subcaptionbox*{}{\rotatebox[origin=t]{90}{\small \hspace{.5in} }}
  \subcaptionbox{\small $K \in \{2,\, 20,\, 200\}$, $\rho = .02$}
{\includegraphics[width=.2\textwidth,angle=-90]{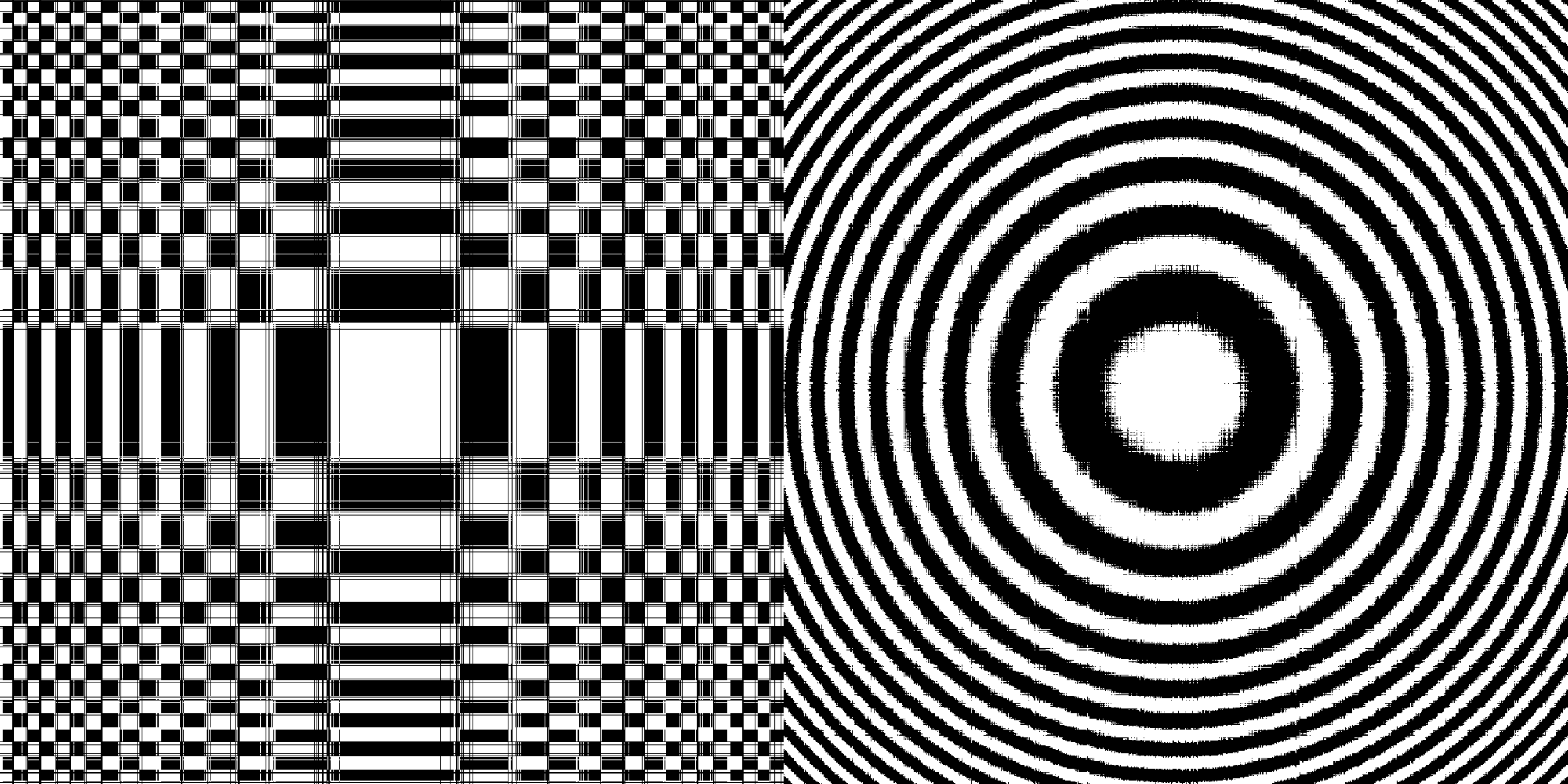}
%  \subcaptionbox{\small $K = 20$}
\includegraphics[width=.2\textwidth,angle=-90]{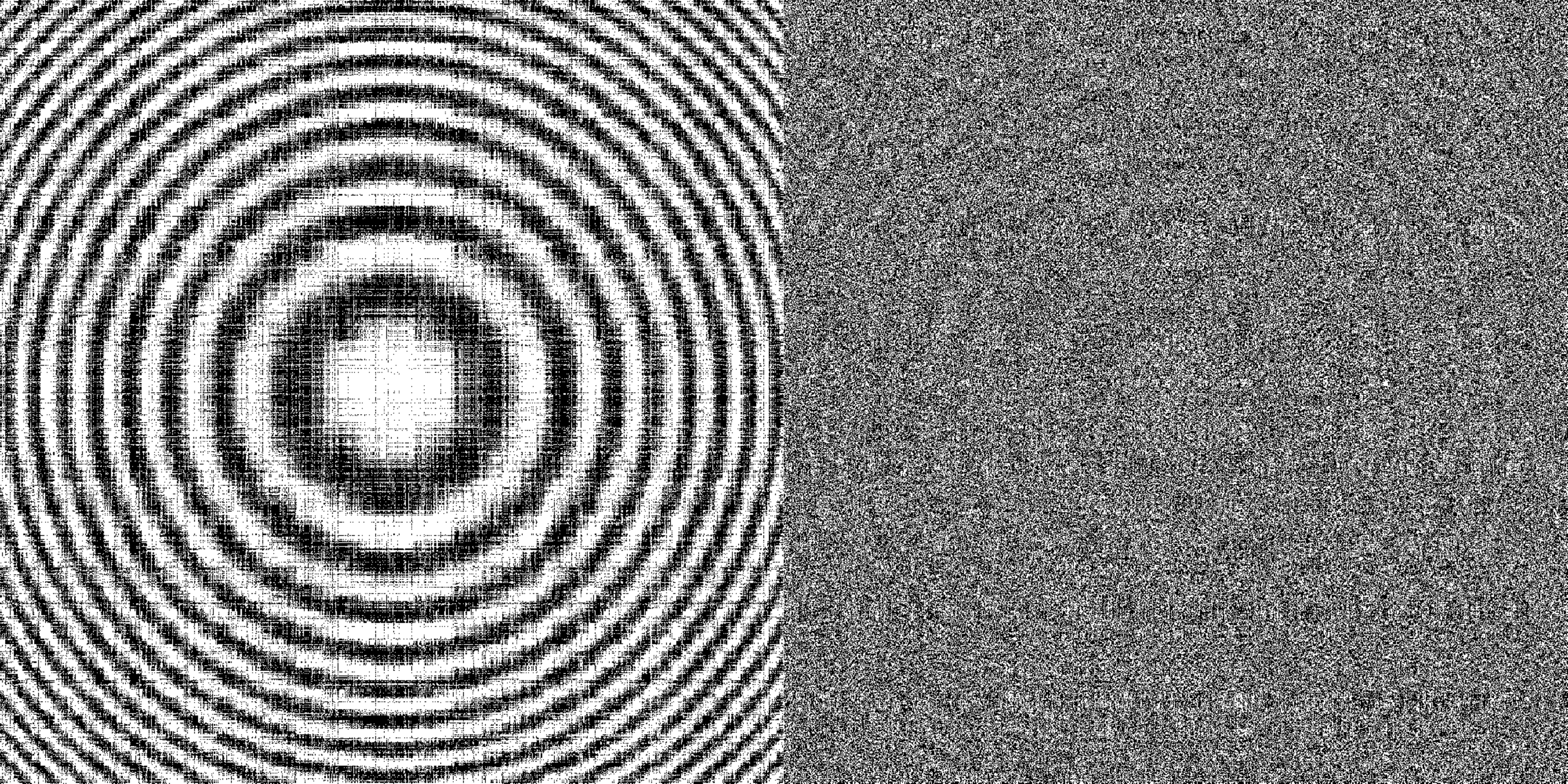}
%  \subcaptionbox{\small $K = 200$}
\includegraphics[width=.2\textwidth,angle=-90]{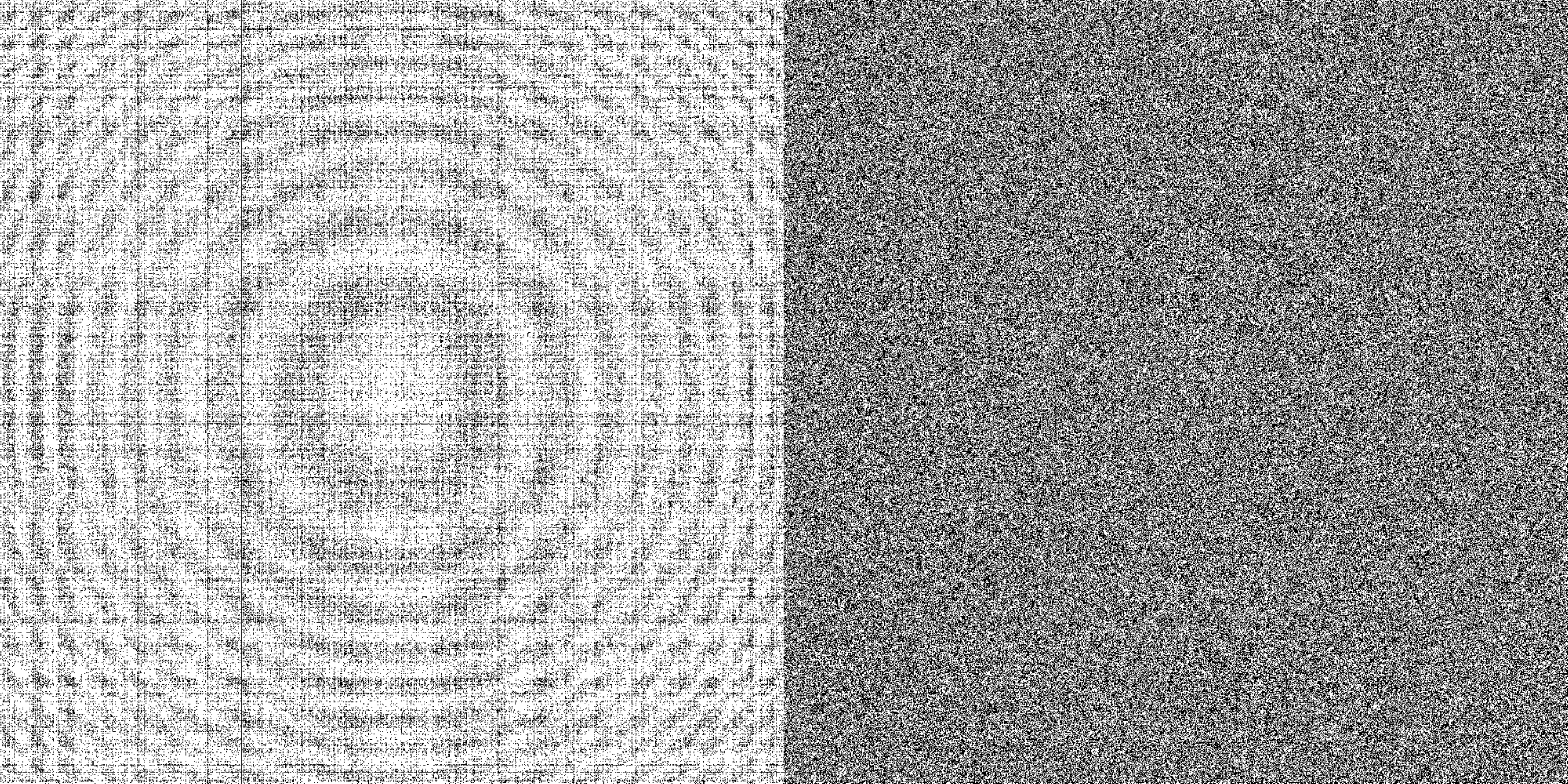}}
\hspace{.05in}
\subcaptionbox{reg.}{\includegraphics[width=.2\textwidth,angle=-90]{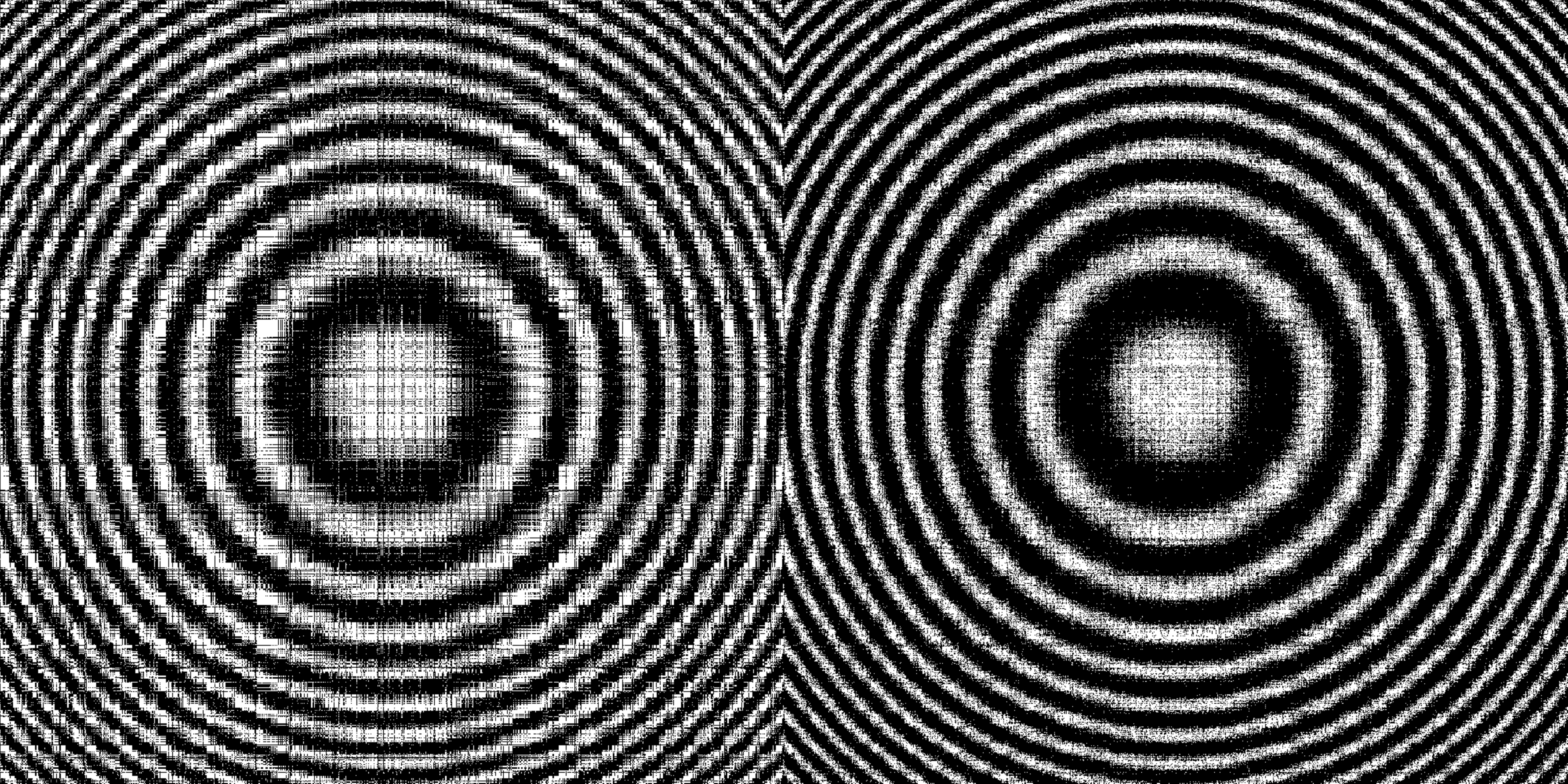}}
}
\vspace{-.1in}
\caption{\it \small  Comparison of low-rank Boolean matrix completion using 1) message passing (using Boolean factors) and 2) GLRM (using real-valued factors) for $\K = 10$. The first column shows the original image (top) and the observation for $\rho = .01$ (bottom).  
\textbf{(a)}~increasing numbers of observations $\rho$; \textbf{(b)}~increasing rank $\K$; 
\textbf{(c)}~using quadratic regularization for 
GLRM and sparsity inducing priors $\PPX_{\m,\k}(0) = \PPX_{\m,\k}(0) = .9$ for message passing. Here $\K = 20$ and $\rho = .02$ -- \ie similar to the figure (b) middle.}
\label{fig:circles}
\end{figure*}

\section{Uninfluential Edges}
\begin{figure}
  \begin{center}
    \includegraphics[width=0.45\textwidth]{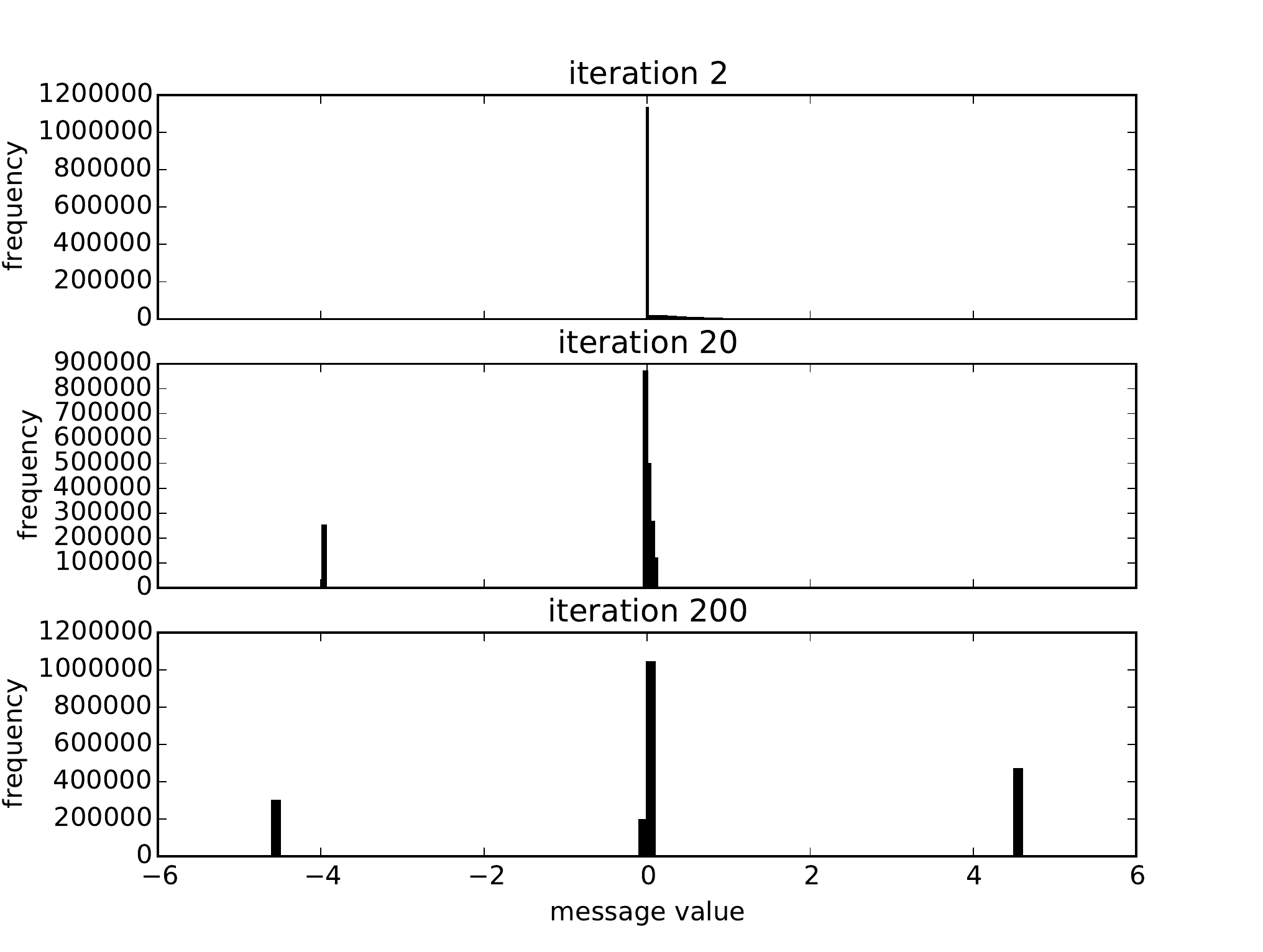}
  \end{center}
  \caption{{\it \small Histogram of BP messages $\{\aAh_{\m,\n}\ttt{t}\}_{\m,\n}$ at $t \in \{2,20,200\}$ for a random $1000 \times 1000$ matrix factorization with $\K = 2$.}}\label{fig:histogram}
\vspace*{-1em}
\end{figure}

\Cref{fig:histogram} shows the histogram of factor-to-variable messages $\{\aAh_{\m,\n}\}_{1\leq \m \M,1\leq \n \leq N}$ at different iterations. 
It suggests that a large portion of messages are close to zero. Since these are log-ratios, the corresponding probabilities are close to uniform. 
Uniform message over an edge in a factor-graph is equivalent to non-existing edges, which 
in turn reduces the number of influential loops in the factor-graph.

\section{Image Completion}\label{sec:image}
\Cref{fig:circles} is an example 
of completing a $1000 \times 1000$ black and white image, here using message passing or GLRM.
In \Cref{fig:circles}\textbf{(a)} we vary the number of observed pixels $\rho \in \{.01,\, .02,\, .05\}$ with fixed $\K = 10$ and 
in \Cref{fig:circles}\textbf{(b)} we vary the rank $\K \in \{2,\,20,\,200\}$, while fixing $\rho = .02$. 
A visual inspection of reconstructions suggests that, since GLRM is using real factors,
it can easily over-fit the observation as we increase the rank.
However, the Boolean factorization, despite being expressive, does not show over-fitting behavior for larger rank values -- as if the result was regularized.
In \Cref{fig:circles}\textbf{(c)}, we regularize both methods for $\K = 20$: 
for GLRM we use Gaussian priors over both $\XX$ and $\YY$ and for message passing
we use sparsity inducing priors $\PPX_{\m,\k}(0) = \PPX_{\m,\k}(0) = .9$. 
This improves the performance of both methods. 
However, note % we would like to emphasize 
that regularization does not significantly improve the results of GLRM when applied to the matrix completion task, where the underlying factors are known to be Boolean (see \Cref{fig:stuff}(right)).

%Here, having twice as much noise corresponds to adding $\log(2)$ to the value of $c$ in 
%$\oR _{\m, \n} = c (\oO_{\m,n} - .5)$}

%%% Local Variables:
%%% mode: latex
%%% TeX-master: "sample_paper"
%%% End:

\end{document}